\documentclass[preprint,12pt]{elsarticle}
\usepackage{multirow}
\usepackage{float}
\usepackage{makecell}
\usepackage{graphicx}
\usepackage{subfigure}
\usepackage{amssymb}
\usepackage{mathtools}
\usepackage{url}
\usepackage{textcomp,gensymb}
\usepackage{algorithm,algpseudocode} 
\usepackage{tikz}
\usepackage{xcolor}
\usepackage{booktabs}
\usepackage{amsmath}
\usepackage{amssymb}
\usepackage{hyperref}

\usepackage{amsthm}

\newtheorem{thm}{Theorem}
\newtheorem{rem}{Remark}

\newtheorem{lemma}{Lemma}

\theoremstyle{definition}
\newtheorem{exam}{Example}

\global\long\def\norm#1{\lVert#1\rVert}%
\global\long\def\trace#1{\text{trace}(#1)}%
\global\long\def\abs#1{\lvert#1\rvert}%
\global\long\def\mtl{\mathcal{L}}%
\global\long\def\mtf{\mathcal{F}}%
\global\long\def\mtp{\mathcal{P}}%
\global\long\def\mtc{\mathcal{C}}%
\global\long\def\mte{\mathcal{E}}%
\global\long\def\esssup{\text{esssup}}%
\global\long\def\mtq{\mathcal{Q}}%

\usepackage{graphicx}
\usepackage{textcomp}

\global\long\def\norm#1{\lVert#1\rVert}%
\global\long\def\abs#1{\lvert#1\rvert}%
\global\long\def\mtf{\mathcal{F}}%
\global\long\def\mtp{\mathcal{P}}%
\global\long\def\mte{\mathcal{E}}%
\begin{document}
\begin{frontmatter}

\title{\LARGE \bf Funnel Synthesis via
LMI Copositivity Conditions for Nonlinear Systems}
\author[label1]{Taewan Kim}
\author[label1]{Beh\c{c}et A\c{c}\i kme\c{s}e}
\affiliation[label1]{organization={Department of Aeronautics and Astronautics, University of Washington},
            addressline={3940 Benton Ln NE}, 
            city={Seattle},
            postcode={98195}, 
            state={WA},
            country={United States}}
            


\begin{abstract}
Funnel synthesis refers to a procedure for synthesizing a time-varying controlled invariant set and an associated control law around a nominal trajectory. The computation of the funnel involves solving a continuous-time differential equation or inequality, ensuring the invariance of the funnel. Previous  approaches often compromise the invariance property of the funnel; for example, they may enforce the equation or the inequality only at discrete temporal nodes and do not have a formal guarantee of invariance at all times. This paper proposes a computational funnel synthesis method that can satisfy the invariance of the funnel without such compromises. We derive a finite number of linear matrix inequalities (LMIs) that imply the satisfaction of  a continuous-time differential linear matrix inequality guaranteeing the invariance of the funnel at all times from the initial to the final time. To this end, we utilize LMI conditions ensuring  matrix copositivity, which then imply continuous-time invariance. The primary contribution of the paper is to prove that the resulting funnel is indeed invariant over a finite time horizon. We validate the proposed method via a three-dimensional trajectory planning and control problem with obstacle avoidance constraints, and a six-degree-of-freedom powered descent guidance.
\end{abstract}

\end{frontmatter}

\section{Introduction}

Funnels, also referred to as tubes in model predictive control (MPC) literature \citep{mayne2005robust,rakovic2016elastic}, represent a time-varying controlled
invariant set around a nominal trajectory over a finite-time horizon.
Funnel synthesis refers to a procedure not only computing the invariant
funnel but also synthesizing an associated feedback control. The
computed funnel with the associated controller can serve various purposes
including robust motion planning
and control \citep{majumdar2017funnel} and trajectory
planning with guaranteed feasibility \citep{reynolds2021funnel}. 

To derive the invariance conditions of the funnel, Lyapunov theory, dissipative theory, or Hamilton-Jacobi (HJ) reachability analysis have been commonly employed
\citep{seo2021fast,reynolds2021funnel}. These approaches result in
differential equations or inequalities that ensure the invariance of the funnel with respect to variables such as ellipsoidal parameters and feedback gains. 
For instance, when employing
a quadratic Lyapunov (storage) function, the invariance condition derived by
Lyapunov theory has the form of a differential linear matrix inequality
(DLMI) \citep{10167750}. 
One challenge is solving the resulting differential equation or DLMI over a finite horizon especially when state and input constraints are involved.

Various methods have been employed to compute the finite-horizon controlled
invariant funnel for nonlinear systems \citep{tobenkin2011invariant,seo2021fast}.
One notable approach is the use of sum-of-squares (SOS) programming \citep{majumdar2017funnel}.
This method has been applied to polynomial Lyapunov functions based
on Lyapunov theory. One issue in the work is that the resulting invariance condition derived by Lyapunov theory is only enforced at discrete node points. The work in \citep{jang2021fast} extends this by ensuring continuous-time invariance by employing the linear interpolation in time and allowing discontinuities at each node point. They compute the forward reachable set of the system for an admissible input set sequentially, within each time interval between discrete node points. However, their approach is not applicable to our problem, where both the feedback control and the controlled-invariant funnel need to be computed in the presence of constraints. In \citep{majumdar2017funnel,jang2021fast}, to apply SOS programming, the system dynamics must be approximated as polynomials, yet these works neglect the error arising from such polynomial approximations and do not provide frameworks to account for this error. More importantly, these studies overlook state and input constraints except for the input saturation, and hence the objective of funnel computation merely focuses on minimizing the size to avoid constraint violations.

In \citep{reynolds2021funnel}, the computation of a quadratic funnel using a quadratic Lyapunov function is studied, directly imposing the state and input constraints.  This approach transforms the resulting DLMI into a finite number of linear matrix inequalities (LMIs), assuming Jacobian matrices to be within known polytopes. However, this assumption generally does not hold when the Jacobian matrices are evaluated around arbitrary dynamically-feasible nominal trajectories. Our previous work \citep{10167750} takes a different approach, computing the quadratic funnel through numerical optimal control using a multiple shooting method, without making the assumption on the Jacobian matrices.  The downside of the previous work is that the resulting DLMI is only ensured at each temporal node, thereby risking the invariance condition violations between these nodes. 

The proposed work synthesizes funnels based on uncertain linear time-varying (LTV) systems obtained via linearization of the original nonlinear dynamics around given nominal trajectories. Hence,  the following papers on the analysis and synthesis of uncertain LTV systems are relevant to this paper: \cite{buch2021finite,seiler2019finite,seiler2024trajectory,pfifer2015robustness}.  The results in \citep{seiler2019finite} provide a computational approach to analyze robust performance of uncertain LTV systems whose uncertainties are characterized by integral quadratic constraints (IQCs). This work is later extended in \cite{buch2021finite} to the synthesis of robust output-feedback controllers. The algorithm proposed in \cite{buch2021finite} alternately computes the variables describing the Lyapunov function for the closed-loop LTV system with a fixed feedback controller. Then it computes the feedback control gains while fixing  the variables defining the underlying Lyapunov function until convergence. That is,  they fix a set of solution variables while solving the rest, and then they solve for the solution variables that they fixed while the rest are now kept constant, and the process repeats in an iterative fashion.  On the other hand, our proposed work characterizes nonlinearities using incremental quadratic constraints based on Lipschitz nonlinearities and constructs a single convex optimization problem whose solution provides both the funnel and feedback controller  parameters.

This paper provides a computational funnel synthesis method for nonlinear systems under bounded disturbances over a finite-time horizon. The main contribution is proving that the computed funnel is indeed invariant over the entire time interval. To achieve this, we first derive an incremental dynamical system representing the behavior relative to the given nominal trajectory,  modeled as an uncertain LTV system. The key idea is to approximate this LTV system to an uncertain linear parameter-varying (LPV) system where the approximation error is treated as a state-, input-, and disturbance-dependent uncertainty. This enables the derivation of a finite number of LMIs whose satisfaction guarantees the required continuous-time DLMI for the invariance of the funnel. In this derivation, we employ LMI conditions that ensure  matrix copositivity \citep{parrilo2000semidefinite,arceo2020copositive}. Including linear
state and input constraints together, the funnel synthesis problem
can be finally formulated as a finite-dimensional semidefinite program
(SDP). Since both the system's nonlinearity and LPV approximation errors are explicitly accounted for, the resulting funnel is guaranteed to be invariant for the original nonlinear system.

\subsection{Notation}

The set notations $\mathbb{R}$, $\mathbb{R}_{+}$, $\mathbb{Z}$
and $\mathbb{R}^{n}$ are the sets of real, nonnegative real, integer,
and the $n$-dimensional Euclidean spaces, respectively. Intervals
are written by $\mathbb{Z}_{[a,b)}=\{z\in\mathbb{Z}:a\leq z<b\}$.
The symmetric matrix $Q=Q^{\top}(\succeq)\succ0$ implies $Q$ is
a positive (semi) definite PD (PSD) matrix. The set of positive (semi) definite
matrices whose size is $n\times n$ are denoted by $\mathbb{S}_{++}^{n}(\mathbb{S}_{+}^{n})$.
The identity matrix having $n\times n$ size is denoted by $I_{n}$.
The subscript and the time argument will be omitted when it is clear
from the context. The notation $\mtl_{2}[a,b]$ is a set of Lebesgue
measurable functions $x(t)$ defined on an interval $[a,b]\subset\mathbb{R}$
such that $\left(\int_{a}^{b}x(t)^{\top}x(t)\text{d}t\right)^{1/2}<\infty.$
A property is said to hold almost everywhere, or for almost every
$x\in X$ with some set $X$ if the set on which it fails is a Lebesgue
measure zero set. We abbreviate the notations $A^{\top}PA$ and $\left[\begin{array}{cc}
a & b^{\top}\\
b & c
\end{array}\right]$ as $(\star)^{\top}PA$ and $\left[\begin{array}{cc}
a & \star\\
b & c
\end{array}\right]$, respectively. The Minkowski sum is denoted by $\oplus$. The stacking of vectors $x$ and $y$ is denoted by $(x, y) = [x^\top y^\top]^\top$. The operation $\text{diag}(\cdot)$ outputs a diagonal matrix formed from a given vector.

\section{Invariance of funnel} \label{sec:2}
\subsection{Nonlinear systems}

Consider the finite-horizon continuous-time nonlinear systems
in the form of
\begin{equation}
\dot{x}(t)=f(t,x(t),u(t),w(t)),\quad\forall t\in[t_{0},t_{f}],\label{eq:nonlinear_system}
\end{equation}
where $x(t)\in\mathbb{R}^{n_{x}}$ is the state, $u(t)\in\mathbb{R}^{n_{u}}$
is the input, $w(t)\in\mathbb{R}^{n_{w}}$ is the (exogenous) disturbance,
and $t_{0},t_{f}$ are the initial and final time. The system dynamics
$f:\mathbb{R}_{+}\times\mathbb{R}^{n_{x}}\times\mathbb{R}^{n_{u}}\times\mathbb{R}^{n_{w}}\rightarrow\mathbb{R}^{n_{x}}$
are assumed to be continuous in $t$ and continuously
differentiable in $x,\, u$ and $w$. We assume
that the input signal $u(\cdot)\in\mtl_{2}^{n_{u}}[t_{0},t_{f}]$
 and the disturbance signal $w(\cdot)\in\mtl_{2}^{n_{w}}[t_{0},t_{f}]$ are piecewise continuous. We further assume that the disturbance signal $w(\cdot)$  is essentially bounded from above by one,
that is, 
\begin{equation}
\norm{w(\cdot)}_{\infty}\leq1,\label{eq:bound_on_2}
\end{equation}
where $\norm{w(\cdot)}_{\infty}\coloneqq\esssup_{t\in[t_{0},t_{f}]}\norm{w(t)}_{2}$. Notice that the choice of the upper bound (i.e., one) is not restrictive, since if $w(\cdot)$ is bounded above by some constant $w_{max}\in\mathbb{R}_{+}$, we can redefine $w\leftarrow w/w_{max}$ so that the bound becomes one. 
We refer to a trajectory as a collection of the state, the input,
and the disturbance signals, denoted together by $(x(\cdot),u(\cdot),w(\cdot))$.
The nominal trajectory having the zero disturbance $(\bar{x}(\cdot),\bar{u}(\cdot),0)$
is assumed to be given and dynamically feasible, that is, $\dot{\bar{x}}=f(t,\bar{x},\bar{u},0)$ for all $t\in[t_{0},t_{f}]$.
With the given nominal trajectory, we can rewrite the system \eqref{eq:nonlinear_system}
in the linear fractional form 
\begin{subequations}
\label{eq:LFT}
\begin{align}
f(t,x,u,w) & =A(t)x+B(t)u+F(t)w+E_{o}\phi(t,q)\\
q & =C_{o}x+D_{o}u+G_{o}w,
\end{align}
\end{subequations}
where the function $\phi:\mathbb{R}_{+}\times\mathbb{R}^{n_{q}}\rightarrow\mathbb{R}^{n_{\phi}}$
represents the nonlinearity of the system and $q(t)\in\mathbb{R}^{n_{q}}$
is its argument. The matrices $A(t)$, $B(t)$, and $F(t)$ are chosen
as first-order approximations of \eqref{eq:nonlinear_system} around
the nominal trajectory. We assume that $A(t)$, $B(t)$, and $F(t)$ are bounded
for all $t\in[t_{0},t_{f}]$. The constant matrices $E_{o}\in\mathbb{R}^{n_{x}\times n_{\phi}}$,
$C_{o}\in\mathbb{R}^{n_{q}\times n_{x}}$, $D_{o}\in\mathbb{R}^{n_{q}\times n_{u}}$,
and $G_{o}\in\mathbb{R}^{n_{q}\times n_{w}}$ are selector matrices, composed of 0s and 1s, chosen to structure the nonlinearity of the
system. 

\begin{exam}\label{exm:unicycle}
The unicycle model can be written as
\begin{equation}
\dot{x}=\left[\begin{array}{c}
\dot{x}_{1}\\
\dot{x}_{2}\\
\dot{x}_{3}
\end{array}\right]=\left[\begin{array}{c}
u_{1}\cos(x_{3}+c_{1}w_{1})\\
u_{1}\sin(x_{3}+c_{1}w_{1})\\
u_{2}+c_{2}w_{2}
\end{array}\right],\label{eq:unicycle}
\end{equation}
where $x_{1}$ and $x_{2}$ are $x$- and $y$-positions, $x_{3}$
are the yaw angle, $u_{1}$ is the velocity control, and $u_{2}$ is the
angular velocity control. The disturbances $w_{1}$ and $w_{2}$ affect
to the yaw angle estimation and the angular velocity control, and
$c_{1},c_{2}\in\mathbb{R}_{+}$ are system parameters. The argument
$q$ of $\phi$ can be chosen as $[x_3,u_1,w_1]$ with
\begin{align*}
C_{o}&=\left[\begin{array}{ccc}
0 & 0 & 1\\
0 & 0 & 0\\
0 & 0 & 0
\end{array}\right], D_{o}=\left[\begin{array}{cc}
0 & 0\\
1 & 0\\
0 & 0
\end{array}\right], G_{o}=\left[\begin{array}{cc}
0 & 0\\
0 & 0\\
1 & 0
\end{array}\right].
\end{align*}
Since the only first two components in $f$ involve with the nonlinearity, the matrix $E_{o}$ is  $\left[\begin{array}{ccc}
1 & 0 & 0\\
0 & 1 & 0
\end{array}\right]^{\top}$.
The time-varying matrices $A(t),\,B(t)$, and $F(t)$ are given as follows:
\begin{subequations}
\label{eq:ABF_unicycle}
\begin{align}
A(t) &=\left[\begin{array}{ccc}
0 & 0 & -\bar{u}_{1}(t)\sin\bar{x}_{3}(t)\\
0 & 0 & \bar{u}_{1}(t)\cos\bar{x}_{3}(t)\\
0 & 0 & 0
\end{array}\right], B(t)=\left[\begin{array}{cc}
\cos\bar{x}_{3}(t) & 0\\
\sin\bar{x}_{3}(t) & 0\\
0 & 1
\end{array}\right], \\
F(t) & =\left[\begin{array}{cc}
-c_{1}\bar{u}_{1}(t)\sin\bar{x}_{3}(t) & 0\\
c_{1}\bar{u}_{1}(t)\cos\bar{x}_{3}(t) & 0\\
0 & c_{2}
\end{array}\right]
\end{align}
\end{subequations}
where $\bar{x}=[\bar{x}_1,\bar{x}_2,\bar{x}_3]^\top$ and $\bar{u}=[\bar{u}_1,\bar{u}_2]^\top$ are the nominal state and input, respectively. The nonlinear function $\phi$ that is the remainder term in \eqref{eq:LFT} is given by 
\begin{align*}
    \phi(t,q) = \left[\begin{array}{c}
c_{1}\bar{u}_{1}w_{1}\sin\bar{x}_{3}-u_{1}\cos\bar{x}_{3}+u_{1}\cos(c_{1}w_{1}+x_{3})+\bar{u}_{1}x_{3}\sin\bar{x}_{3}, \\
-c_{1}\bar{u}_{1}w_{1}\cos\bar{x}_{3}-u_{1}\sin\bar{x}_{3}+u_{1}\sin(c_{1}w_{1}+x_{3})-\bar{u}_{1}x_{3}\cos\bar{x}_{3}
\end{array}\right],
\end{align*}
where the time argument is omitted.
\end{exam}

\subsection{Incremental dynamical system and its LPV approximation}

The incremental form of dynamics \citep{d2013incremental,xu2020observer,accikmecse2011robust} illustrates the behavior of the system
\eqref{eq:nonlinear_system} relative to the nominal trajectory.
To derive it, we first define difference variables as
\begin{align*}
\eta(t)\coloneqq x(t)-\bar{x}(t),\quad\xi(t)=u(t)-\bar{u}(t),\quad\delta q(t)=q(t)-\bar{q}(t),    
\end{align*}
where $\bar{q}(t)\coloneqq C_{o}\bar{x}(t)+D_{o}\bar{u}(t)$. Having $\eta$
as the state, the incremental dynamics can be derived as the following
uncertain LTV system:
\begin{subequations}
\label{eq:LTV_1}
\begin{align}
\dot{\eta} & =A(t)\eta+B(t)\xi+F(t)w+E_{o}\left(\phi(t,q)-\phi(t,\bar{q})\right),\\
\delta q & =C_{o}\eta+D_{o}\xi+G_{o}w.
\end{align}
\end{subequations}
Since the original system $f$ in \eqref{eq:nonlinear_system} is
continuously differentiable, $f$ is locally Lipschitz, thereby resulting
in $\phi$ being locally Lipschitz with its second argument. Then, for each $t\in[t_{0},t_{f}]$,
there exists a local Lipschitz constant $\gamma_i(t)\in\mathbb{R}_{+}$
such that
\begin{align}
\norm{\phi_i(t,q)-\phi_i(t,\bar{q})}_{2}\leq\gamma_i(t)\norm{q(t)-\bar{q}(t)}_{2},\quad\forall i \in \mathbb{Z}_{[1,n_\phi]},\quad\forall q,\bar{q}\in\mtq,\label{eq:locally_lipschitz}
\end{align}
for any compact set $\mtq\subseteq \mathbb{R}^{n_q}$. Here the term $\phi(t,q)-\phi(t,\bar{q})$,
due to the nonlinearity of \eqref{eq:LFT}, is a state-, input-, and disturbance-dependent uncertainty characterized
by the constraint \eqref{eq:locally_lipschitz}. This type of LTV system in incremental form with Lipschitz nonlinearity has been studied in the context of funnel computation in several existing works. See, for example, \cite{reynolds2021funnel,10167750} for related formulations and applications.

Now we describe an LPV approximation of the LTV system \eqref{eq:LTV_1}. The underlying motivation of this approximation is that we can derive a finite number of LMIs that imply the invariance of the funnel for the LPV system, which is not tractable for the LTV system. We start by partitioning the time horizon $[t_{0},t_{f}]$ into $N\in\mathbb{R}_+$ uniform subintervals using time nodes defined as:
\[
t_{k}=t_{0}+\frac{k}{N}(t_{f}-t_{0}),\quad\forall k\in\mathbb{Z}_{[0,N]}.
\]
On each subinterval $[t_k, t_{k+1}]$, we approximate the time-varying system matrices $A(t)$, $B(t)$,
and $F(t)$ using a first-order hold (FOH) approach, which linearly interpolates between their values at the endpoints of the interval. This results in convex combinations of the matrices at $t_k$ and $t_{k+1}$. We define:
\begin{subequations}
\label{eq:FOH}
\begin{align}
\tilde{\square}(t) & \coloneqq \sigma_{1}^{k}(t)\square_{k}+\sigma_{2}^{k}(t)\square_{k+1},\quad t\in[t_k,t_{k+1}]\\
\sigma_{1}^{k}(t) & =\frac{t_{k+1}-t}{t_{k+1}-t_{k}},\quad\sigma_{2}^{k}(t)=\frac{t-t_{k}}{t_{k+1}-t_{k}},\label{eq:FOH_sigma}
\end{align}
\end{subequations}
where $\square_{k}=\square(t_{k})$ and the placeholder $\square$ corresponds to $A$, $B$, and $F$. By applying \eqref{eq:FOH} across all $k$ in $\mathbb{Z}_{[0,N-1]}$, we obtain continuous piecewise linear approximations $(\tilde{A}(t),\tilde{B}(t),\tilde{F}(t))$ of the system matrices $(A(t),B(t),F(t))$ over the entire horizon $[t_0,t_f]$. 

It is worth noting that the system matrices $(A(t),B(t),F(t))$ are not necessarily piecewise linear, and therefore the approximation error is generally nonzero.
However, it is common in the literature (e.g., \cite{reynolds2021funnel,malikov2020numerical}) to assume piecewise linearity to simplify funnel computation. In contrast, this paper does not assume zero approximation error; instead, the error is modeled explicitly as a state-, input-, and disturbance-dependent uncertainty.

Using \eqref{eq:FOH},
we can rewrite \eqref{eq:LTV_1} equivalently as follows:
\begin{subequations}
\label{eq:LPV_1}
\begin{align}
\dot{\eta} & =\tilde{A}(t)\eta+\tilde{B}(t)\xi+\tilde{F}(t)w+\underbrace{\Delta_e(t)\left[\begin{array}{c}
\eta\\
\xi\\
w
\end{array}\right]}_{\coloneqq e(t)}+E_{o}\delta\phi,\label{eq:LPV_1_1}\\
\Delta_e(t) & \coloneqq\left[\begin{array}{ccc}
A(t)-\tilde{A}(t) & B(t)-\tilde{B}(t) & F(t)-\tilde{F}(t)\end{array}\right],\label{eq:LPV_1_2}
\end{align}
\end{subequations}
where $\delta\phi(t,\delta q)\coloneqq\phi(t,q)-\phi(t,\bar{q})$.
Note that two representations \eqref{eq:LTV_1} and \eqref{eq:LPV_1} are equivalent by the definition of the error term $e(t)$ in \eqref{eq:LPV_1_1}.
Since $f$ is assumed to be continuously differentiable, the matrices $A(t)$, $B(t)$, and $F(t)$ are continuous, so bounded within the finite interval $[t_0,t_f]$. Hence,
 there exists a positive constant $\beta_e(t)\in\mathbb{R}_{+}$ for
each $t$ such that
\begin{align}
\norm{\Delta_e(t)}_{2}\leq\beta_e(t),\label{eq:bound_on_Delta}
\end{align}
where $\norm{\cdot}_2$ denotes the matrix 2-norm.

Alternatively, we can express the error term $e(t)$ more compactly as
\begin{align*}
    e(t) = E_\Delta \Delta (t) q_{\Delta}(t)
\end{align*}
where the matrix $E_\Delta \in \mathbb{R}^{n_x \times n_\Delta}$, the block-diagonal matrix $\Delta \in \mathbb{R}^{n_\Delta \times n_{q_\Delta}}$  of uncertainties, and the vector $q_\Delta \in \mathbb{R}^{n_{q_\Delta}}$ are given by
\begin{align*}
    E_\Delta &= \left[\begin{array}{ccc} E_\eta & E_\xi & E_w \end{array}\right],\quad \Delta = \left[\begin{array}{ccc} 
    \Delta_\eta & 0 & 0 \\
    0 & \Delta_\xi & 0 \\
    0 & 0 & \Delta_w
    \end{array}\right], \quad
    q_\Delta = \left[\begin{array}{c} q_\eta \\ q_\xi \\ q_w \end{array}\right],
\end{align*}
 with $n_\Delta \in \mathbb{R}_{+}$ and $n_{q_\Delta} = n_{q_\eta}+n_{q_\xi}+n_{q_w}$. These components are chosen such that
\begin{align*}
    (A(t)-\tilde{A}(t))\eta &= E_\eta \Delta_\eta q_\eta, \\
    (B(t)-\tilde{B}(t))\xi &= E_\xi \Delta_\xi q_\xi, \\
    (F(t)-\tilde{F}(t))w &= E_w \Delta_w q_w,
\end{align*}
where the vectors $q_\eta$, $q_\xi$, and $q_w$ are linear functions of $\eta$, $\xi$, and $w$, respectively, as follows:
\begin{align*}
    q_\eta = C_\eta \eta,\quad q_\xi = D_\xi \xi,\quad q_w = G_w w.
\end{align*}
In this formulation, the uncertainties are captured using the following inequalities:  
\begin{align}
\norm{\Delta_\square(t)}_2 \leq \beta_\square(t),\label{eq:bound_on_Delta_diag}
\end{align}
where $\beta_\square(t)>0$ bounds the corresponding uncertainty $\Delta_\square(t)$, with the placeholder $\square$ representing the subscripts $\eta$, $\xi$, and $w$. Compared to the single norm-bounded inequality in \eqref{eq:bound_on_Delta}, these bounds in \eqref{eq:bound_on_Delta_diag} can represent the uncertainly more compactly by exploiting its structure. More details on the structured uncertainty representation can be found in \cite[6.2.1]{reynolds2020computational}, \cite{doyle1982analysis}, \cite{boyd1994linear}.
    

Now we can summarize the incremental dynamics, written in the form of an uncertain LPV system, as follows:
\begin{subequations}
\label{eq:LPV_final}
\begin{align}
\dot{\eta} & =\tilde{A}(t)\eta+\tilde{B}(t)\xi+\tilde{F}(t)w+Ep(t,r), \label{eq:LPV_system}\\
r & \coloneqq C\eta+D\xi+Gw=[q_\Delta^\top,\delta q^{\top}]^\top\\
p & \coloneqq\left[\begin{array}{c}
\Delta(t) q_\Delta (t) \\
\delta\phi(t,\delta q)
\end{array}\right],\quad E\coloneqq\left[\begin{array}{cc}
E_\Delta & E_{o}\end{array}\right],\\
C & \coloneqq\left[\begin{array}{c}
C_\eta\\
0\\
0\\
C_{o}
\end{array}\right], D\coloneqq\left[\begin{array}{c}
0\\
D_\xi\\
0\\
D_{o}
\end{array}\right], G\coloneqq\left[\begin{array}{c}
0\\
0\\
G_w\\
G_{o}
\end{array}\right].
\end{align}
\end{subequations}
 All uncertainties in the system \eqref{eq:LPV_final} are lumped
into the vector $p(t,r)\in\mathbb{R}^{n_{p}}$ where $n_{p}=n_{\Delta}+n_{\phi}$
and $r(t)\in\mathbb{R}^{n_{r}}$. We characterize the uncertain term
$p$ using the following quadratic inequality (QI) \citep{accikmecse2003robust,accikmecse2011robust} for each $t$:
\begin{equation}
\left[\begin{array}{c}
r(t)\\
p(t,r)
\end{array}\right]^{\top}M(t)\left[\begin{array}{c}
r(t)\\
p(t,r)
\end{array}\right]\geq0,\label{eq:quadratic_inequality}
\end{equation}
where the matrix $M(t)\in\mathbb{R}^{(n_{r}+n_{p})\times(n_{r}+n_{p})}$
is called a multiplier matrix. It follows from \eqref{eq:locally_lipschitz} and \eqref{eq:bound_on_Delta} that the valid multiplier matrix has the form of
\begin{subequations}
\label{eq:multiplier_matrix}
\begin{align}
M(t) & =\left[\begin{array}{cc}
N_{1}(t)^{-1} & 0\\
0 & -N_{2}(t)^{-1}
\end{array}\right],
N_{1}(t) \coloneqq\left[\begin{array}{cc}
N_1^\beta (t) & 0\\
0 & \frac{1}{n_\phi} \lambda_{\gamma}(t)I_{n_{q}}
\end{array}\right],\\ N_{2}(t)&\coloneqq\left[\begin{array}{cc}
N_2^\beta (t) & 0\\
0 & \lambda_{\gamma}(t)\Gamma(t)^{2}
\end{array}\right],
N_1^\beta (t) \coloneqq \text{diag}(\lambda_{\beta_\eta}(t) I,\lambda_{\beta_\xi} (t)I,\lambda_{\beta_w}(t) I), \\
N_2^\beta (t) &\coloneqq \text{diag}(\lambda_{\beta_\eta} (t)\beta_\eta^2(t) I,\lambda_{\beta_\xi} (t) \beta_\xi^2 (t) I,\lambda_{\beta_w}(t) \beta_w^2(t) I),
\end{align}
\end{subequations}
where $\Gamma(t) = \text{diag}(\gamma_1(t),\ldots,\gamma_{n_\phi}(t))$ and any positive real-valued functions $\lambda_{\beta_\square}(t)>0$ with the placeholder $\square$ representing $\eta$, $\xi$, and $w$,
and $\lambda_{\gamma}(t)>0$. 
This can be verified by deriving the following:
\begin{align*}
\left[\begin{array}{c}
r\\
p
\end{array}\right]^{\top}M\left[\begin{array}{c}
r\\
p
\end{array}\right] & = \sum_{\square=\{\eta,\xi,w\}}
\frac{1}{\lambda_{\beta_\square}} (q_\square^\top q_\square - \frac{1}{\beta_\square^2} q_\square^\top \Delta_\square^\top \Delta_\square q_\square)\\&+\frac{1}{\lambda_{\gamma}}\left(n_\phi \delta q^{\top}\delta q-\delta\phi^{\top}\Gamma(t)^{-2}\delta\phi\right) \geq 0,
\end{align*}
where the last inequality comes from \eqref{eq:locally_lipschitz} and \eqref{eq:bound_on_Delta}.
The resulting interconnection in the uncertain LPV system is illustrated in Fig.~\ref{fig:interconnection}.

\begin{figure}
    \vspace*{0.3cm}
    \centering
    \begin{tikzpicture}
        \node[draw, rectangle, minimum width=3cm, minimum height=1.5cm, align=center, rounded corners] (box1) at (0,2) {Uncertainties \\ satisfying \eqref{eq:quadratic_inequality} };
        \node[draw, rectangle, minimum width=3cm, minimum height=1.5cm, align=center, rounded corners] (box2) at (0,0) {LPV System \\ \eqref{eq:LPV_system}};
        \node[draw, rectangle, minimum width=3cm, minimum height=1.5cm, align=center, rounded corners] (box3) at (0,-2) {Feedback Control \\ \eqref{eq:linear_feedback}};

        \draw[-latex] (box1.west) -- ++(-1,0) |- ([yshift=0.5cm]box2.west) node[near start, left] {$p$};
        \draw[-latex] ([yshift=0.5cm]box2.east) -- ++(1,0) |- (box1.east) node[near start, right] {$r$};
        
        \draw[-latex] (-3,0) -- (box2.west) node[midway, left, xshift=-1.0cm, align=center] {$w$};
        \draw[-latex] (box2.east) -- (3,0) node[midway, right, xshift=1.0cm, align=center] {$\eta$};

        \draw[-latex] (box3.west) -- ++(-1,0) |- ([yshift=-0.5cm]box2.west) node[near start, left] {$\xi$};
        \draw[-latex] ([yshift=-0.5cm]box2.east) -- ++(1,0) |- (box3.east) node[near start, right] {$\eta$};
    \end{tikzpicture}
    \caption{Uncertain LPV system interconnection with feedback control.}
    \label{fig:interconnection}
\end{figure}

\subsection{Invariance condition}
We consider a scalar-valued continuous Lyapunov function $V:\mathbb{R}\times\mathbb{R}^{n_{x}}\rightarrow\mathbb{R}$
defined by
\begin{equation}
V(t,\eta(t))=\eta(t)^{\top}Q(t)^{-1}\eta(t),\label{eq:Lyapunov_function}
\end{equation}
where $Q(t)\in\mathbb{S}_{++}^{n_{x}}$ is a continuous-time PD matrix-valued continuous function. A \emph{state funnel} is defined as a 1-sublevel
set of $V$, that is,
\begin{equation}
\mte_{Q}(t)\coloneqq\{\eta\mid\eta^{\top}Q(t)^{-1}\eta\leq1\}.\label{eq:state_funnel}
\end{equation}
We employ the following linear time-varying feedback controller:
\begin{align}
    \xi(t)=K(t)\eta(t), \label{eq:linear_feedback}
\end{align}
for the system \eqref{eq:LPV_final} with the feedback gain
$K(t)\in\mathbb{R}^{n_{u}\times n_{x}}$. With this linear feedback,
$\xi$ becomes linear in $\eta$. It follows that
$\eta\in\mte_{Q}$ implies $\xi\in\mte_{KQK^{\top}}$ where
\begin{equation}
\mte_{KQK^{\top}}(t)\coloneqq\{(K(t)Q(t)K(t)^\top)^{\frac{1}{2}}y\mid\norm y_{2}\leq1,y\in\mathbb{R}^{n_{u}}\}.\label{eq:input_funnel}
\end{equation}
We refer $\mte_{KQK^{\top}}$ to as an \emph{input funnel}. We use the subscript $KQK^{\top}$ because, when $KQK^{\top}$ is PD, the set $\mte_{KQK^{\top}}$ is
equivalent to $\{\xi\mid\xi^{\top}(KQK^{\top})^{-1}\xi\leq1\}$. When $KQK^\top$ is only PSD (not PD), $\mte_{KQK^{\top}}$ becomes a degenerated ellipsoid, which remains well-defined and compatible with \eqref{eq:input_funnel}.
Using the state and input funnels, a \emph{funnel} $\mtf(t)$ is defined as an Cartesian product of state and input funnel centered around the nominal state $\bar{x}(t)$ and input $\bar{u}(t)$, respectively, that can be written as
\begin{align}
    \mtf(t) \coloneqq \left(\{\bar{x}(t)\}\oplus \mte_Q(t)\right) \times \left(\{\bar{u}(t)\}\oplus \mte_{KQK^\top}(t)\right).\label{eq:deffunnel}
\end{align}

\begin{lemma}
\label{lem:invariance}Suppose there exists a PD matrix-valued continuous
function $Q(t)\in\mathbb{S}_{++}^{n_{x}}$ and a continuous matrix-valued
function $K(t)\in\mathbb{R}^{n_{u}\times n_{x}}$ such that all trajectories of the system \eqref{eq:LPV_final} satisfy 
\begin{align}
\text{If }  V\geq w^{\top}w, \ \text{then } \dot{V}(t,\eta(t))\leq0,\label{eq:Lyapunov_condition}
\end{align}
for almost all $t\in[t_{0},t_{f}]$. Then, with every almost everywhere continuous signal $w(\cdot)$ such
that \eqref{eq:bound_on_2}, $\mathcal{E}_{Q}(t)$ defined in \eqref{eq:state_funnel}
is an invariant time-varying ellipsoid for the system \eqref{eq:LPV_final}
on $[t_{0},t_{f}]$, that is, if $\eta(\cdot)$ is a solution of \eqref{eq:LPV_final}
with $\eta(t_{0})\in\mte_{Q}(t_{0})$, then $\eta(t)\in\mathcal{E}_{Q}(t)$
for all $t\in[t_{0},t_{f}]$. 
\end{lemma}

\begin{proof}
We prove by contradiction. Suppose that the
invariance condition does not hold. Then, there exists a solution
$\eta(\cdot)$ of \eqref{eq:LPV_final} such that $\eta(t_{0})\in\mte_{Q}(t_{0})$
and $\eta(t_{2})\notin\mte_{Q}(t_{2})$ for some $t_{2}>t_{0}$. Since
$V(t,\eta)$ is continuous in $t$ by continuity of $\eta$ and $Q$, there
exists $t_{1}\in(t_{0},t_{2})$ such that $V(t_{1},\eta(t_{1}))=1$
and $V(t,\eta(t))>1$ for all $t\in(t_{1},t_{2}]$ by the intermediate
value theorem. It follows from \eqref{eq:Lyapunov_condition} that
$\dot{V}(t,\eta(t))\leq0$ for almost all $t\in[t_{1},t_{2}]$ as
\eqref{eq:bound_on_2}. Observe that
\begin{align*}
V(t_{2},\eta(t_{2})) & =V(t_{1},\eta(t_{1}))+\int_{t_{1}}^{t_{2}}\dot{V}(t,\eta(t))\text{d}t,\\
 & \leq V(t_{1},\eta(t_{1})) = 1.
\end{align*}
This implies $\eta(t_{2})\in\mte_{Q}(t_{2})$ that contradicts our
hypothesis that $\eta(t_{2})\notin\mte_{Q}(t_{2})$. Thus, $\mte_{Q}(t)$
is invariant.
\end{proof}

For a detailed discussion of similar Lyapunov conditions in the context of input-to-state stability, we refer the reader to \citep{sontag1995characterizations}. The similar Lyapunov condition with the time-invariant Lyapunov function is also studied in \citep[Lemma 1]{accikmecse2002robust}. While Lemma \ref{lem:invariance} follows a similar proof, we extend the result to consider more general types of disturbance signals $w(\cdot)$ that are continuous almost everywhere.
In the next lemma, we define the invariance of the funnel $\mtf$ and its relevance to the invariance of $\mte_Q$.

\begin{lemma}\label{lem:nonlinear_invariance}
If $\mte_{Q}(t)$ is invariant for the uncertain LPV system \eqref{eq:LPV_final}
with the linear feedback control $\xi=K\eta$, then the 
funnel $\mtf(t)$ defined in \eqref{eq:deffunnel} is invariant for the original
system \eqref{eq:nonlinear_system} in a sense that if $x(\cdot)$ is a solution of \eqref{eq:nonlinear_system} with $x(t_0) \in \{ \bar{x} (t_0) \} \oplus \mte_Q(t_0)$ and the control law,
\begin{equation}
u(t)=\bar{u}(t)+K(t)(x(t)-\bar{x}(t)),\label{eq:control_law}
\end{equation}
then $(x(t),u(t)) \in \mtf(t)$ for all $t\in [t_0,t_f]$.
\end{lemma}

\begin{proof}
Consider a solution $x(\cdot)$ of \eqref{eq:nonlinear_system} with
$x(t_{0})\in\{\bar{x}(t_{0})\}\oplus\mte_{Q}(t_{0})$. This implies
that $\eta(t_{0})\in\mte_{Q}(t_{0})$, and then by the invariance
of $\mte_{Q}$, we have $\eta(t)=x(t)-\bar{x}(t)\in\mte_{Q}(t)$ with
$\xi(t)=u(t)-\bar{u}(t) = K(t)\eta(t)$. Thus, $(x(t),u(t)) \in \mtf(t)$ for all $t\in [t_0,t_f]$
 with \eqref{eq:control_law}.
\end{proof}

\subsection{Differential linear matrix inequality}

This subsection describes the derivation of the DLMI that implies the
invariance condition \eqref{eq:Lyapunov_condition}. 
\begin{lemma}
\label{lem:DLMI}Consider the uncertain LPV system in \eqref{eq:LPV_final}
with the multiplier matrix $M(t)$ satisfying \eqref{eq:multiplier_matrix}.
Suppose there exist $Q(t)\in\mathbb{S}_{++}^{n_{x}}$, $Y(t)\in\mathbb{R}^{n_{u}\times n_{x}}$,
$\lambda_{\beta}(t)\in\mathbb{R}_{+}$, $\lambda_{\gamma}(t)\in\mathbb{R}_{+}$,
and $\lambda_{w}\in\mathbb{R}_{+}$ such that the following differential
matrix inequality holds for some $t\in[t_0,t_f]$:
\begin{align}
D_{LMI}(t,\dot{Q},Q,Y,\lambda_{\beta},\lambda_{\gamma})  &\coloneqq \nonumber\\
\left[\begin{array}{cccc}
W-\dot{Q} & \star & \star & \star\\
N_{2}E^{\top} & -N_{2} & \star & \star\\
\tilde{F}^{\top} & 0 & -\lambda_{w}I & \star\\
CQ+DY & 0 & G & -N_{1}
\end{array}\right]&\preceq0, \label{eq:DLMI}
\end{align}
where $W=\tilde{A}Q+\tilde{B}Y+Q\tilde{A}^{\top}+Y^{\top}\tilde{B}^{\top}+\lambda_{w}Q$ and $N_1,\,N_2$ are defined in \eqref{eq:multiplier_matrix}.
Let $K(t)=Y(t)Q(t)^{-1}$, and then the condition \eqref{eq:Lyapunov_condition}
holds for the given $t$.
\end{lemma}

\begin{proof}
Pre- and post-multiplying \eqref{eq:DLMI} by $\text{diag}(Q^{-1},N_{2}^{-1},I,I)$
generates
\[
\left[\begin{array}{cccc}
Q^{-1}(W-\dot{Q})Q^{-1} & \star & \star & \star\\
E^{\top}Q^{-1} & -N_{2}^{-1} & \star & \star\\
\tilde{F}^{\top}Q^{-1} & 0 & -\lambda_{w}I & \star\\
C_{cl} & 0 & G & -N_{1}
\end{array}\right]\preceq0,
\]
where $C_{cl}=C+DK$. By applying Schur complement,
we can derive
\begin{align*}
\left[\begin{array}{ccc}
Q^{-1}(W-\dot{Q})Q^{-1} & \star & \star\\
E^{\top}Q^{-1} & 0 & \star\\
\tilde{F}^{\top}Q^{-1} & 0 & -\lambda_{w}I
\end{array}\right] +(\star)^{\top}\left[\begin{array}{cc}
N_{1}^{-1} & 0\\
0 & N_{2}^{-1}
\end{array}\right]\left[\begin{array}{ccc}
C_{cl} & 0 & G\\
0 & I & 0
\end{array}\right]\preceq0.
\end{align*}
By pre- and post-multiplying by $[\eta^\top, p^\top, w^\top]$ and its transpose, respectively, and applying the identity $\dot{Q}^{-1} = -Q^{-1} \dot{Q} Q^{-1}$, we obtain
\begin{align*}
\dot{V}+(\star)^{\top}\left[\begin{array}{cc}
N_{1}^{-1} & 0\\
0 & N_{2}^{-1}
\end{array}\right]\left[\begin{array}{c}
r\\
p
\end{array}\right]+\lambda_{w}(V-w^{\top}w)\leq0,
\end{align*}
for all $\eta\in\mathbb{R}^{n_{x}}$, $p\in\mathbb{R}^{n_{p}}$, and
$w\in\mathbb{R}^{n_{w}}$. By using \eqref{eq:quadratic_inequality}, S-procedure \citep{boyd1994linear}, and $\lambda_w > 0$, we can conclude that the above
inequality implies \eqref{eq:Lyapunov_condition}.
\end{proof}

We omit the time arguments in Lemma \ref{lem:DLMI} and its proof for the notational brevity. The matrices $C,\,D,\,E,\,G$, and the constant $\lambda_w$ are time-invariant, while the terms $\tilde{A},\,\tilde{B},\,\tilde{F},\,Q,\,Y,\,K,\,W,\,N_1$, and $N_2$ are time-varying. Notice that with the fixed positive constant $\lambda_{w}$, $D_{LMI}(\dot{Q},Q,Y,\lambda_{\beta},\lambda_{\gamma})$ in \eqref{eq:DLMI}
is linear in its arguments.

\section{Funnel synthesis problem}\label{sec:3}
\subsection{Transformation of DLMI into a finite number of LMIs}

To transform the DLMI \eqref{eq:DLMI} into a finite number of LMIs,
we first parameterize $Q(t)$, $Y(t)$, $\lambda_{\beta}(t)$, and $\lambda_{\gamma}(t)$ by employing the same FOH interpolation used in \eqref{eq:FOH}, that is,
\begin{align*}
\tilde{\square}(t) & \coloneqq \sigma_{1}^{k}(t)\square_{k}+\sigma_{2}^{k}(t)\square_{k+1},\quad t\in[t_k,t_{k+1}],
\end{align*}
where a placeholder $\square$ corresponds to $Q,\,Y,\,\lambda_\beta$, and $\lambda_\gamma$, and the time-varying parameters $\sigma_1^k(t)$ and $\sigma_2^k(t)$ are given in \eqref{eq:FOH_sigma}.  Here
$Q_{k}\in\mathbb{S}_{++}^{n_{x}},Y_{k}\in\mathbb{R}^{n_{u}\times n_{x}}$,
$(\lambda_\beta)_{k}\in\mathbb{R}_{+}$, and $(\lambda_\gamma)_{k}\in\mathbb{R}_{+}$
for all $k\in\mathbb{Z}_{[0,N]}$ are our decision variables to be optimized. With this interpolation, the derivative of $\dot{Q}$ for each open subinterval
$(t_{k},t_{k+1})$, denoted by $\dot{Q}^{k}$, is constant and given by
\begin{align}
\dot{Q}^{k}\coloneqq\dot{Q}=\frac{Q_{k+1}-Q_{k}}{t_{k+1}-t_{k}},\label{eq:FOH_Qdot}
\end{align}
where the superscript $k$ indicates that $\dot{Q}$ is specific to each subinterval.


To derive a finite number of LMIs that ensure the invariance, we also need to represent the time-varying constants $\gamma(t)$ and $\beta(t)$ with a finite set of values. To this end, we define
\begin{align}
    (\gamma_i)_k \coloneqq \sup_{t\in[t_k,t_{k+1}]} \gamma_i(t),\quad \beta_k \coloneqq \sup_{t\in[t_k,t_{k+1}]} \beta(t),
\end{align}
for all $k$ in $\mathbb{Z}_{[0,N-1]}$. These constants $(\gamma_i)_k$ and $\beta_k$ serve as upper bounds for $\gamma_i(t)$ and $\beta(t)$, respectively, on each subinterval $[t_k,t_{k+1}]$. Hence, they become valid constants for \eqref{eq:locally_lipschitz} and \eqref{eq:bound_on_Delta}, respectively on each subinterval $[t_k,t_{k+1}]$. This allows us to incorporate $(\gamma_i)_k$ and $\beta_k$ as constants in the LMIs associated with each subinterval, facilitating the derivation of a finite set of LMIs for the entire horizon $[t_0,t_f]$

With the FOH interpolation \eqref{eq:FOH} for the system matrices and our decision variables, the DLMI condition \eqref{eq:DLMI} for open subinterval $(t_k,t_{k+1})$ can be equivalently rewritten as
\begin{equation}
\sigma_{1}^k(t)\sigma_{1}^k(t)H_{k,k}^k+\sigma_{2}^k(t)\sigma_{2}^k(t)H_{k+1,k+1}^k
+\sigma_{1}^k(t)\sigma_{2}^k(t)(H_{k,k+1}^k+H_{k+1,k}^k)\succeq0,\label{eq:DLMI_parameterization}
\end{equation}
where $H_{i,j}^k\in\mathbb{R}^{n_H\times n_H}$ with $n_H = n_x + n_p + n_w + n_r$ for $i,j\in \{k,k+1\}$ is given by
\begin{align}
H_{i,j} & \coloneqq-\left[\begin{array}{cccc}
W_{i,j}-\dot{Q^k} & \star & \star & \star\\
N_{2,j}^kE^{\top} & -N_{2,j}^k & \star & \star\\
F_{j}^{\top} & 0 & -\lambda_{w}I & \star\\
L_{j} & 0 & G & -N_{1,j}
\end{array}\right],\label{eq:H}\\
W_{i,j} & \coloneqq A_{i}Q_{j}+B_{i}Y_{j}+Q_{j}A_{i}^{\top}+Y_{j}^{\top}B_{i}^{\top}+\lambda_{w}Q_{j},\nonumber\\
L_{j} & \coloneqq CQ_{j}+DY_{j},
N_{1,j}  \coloneqq\text{diag}((\lambda_\beta)_{j}I,(\lambda_\gamma)_{j}I), \nonumber\\
N_{2,j}^k&\coloneqq\text{diag}(\beta_j^{2}(\lambda_\beta)_{j}I,(\lambda_\gamma)_{j}\Gamma_k^{2}),\nonumber
\end{align}
where $\Gamma_k = \text{diag}((\gamma_1)_k,\ldots,(\gamma_{n_\phi})_k)$ and the superscript $k$ is introduced for the piecewise constant terms $\dot{Q}^k$, $\Gamma_k$, and $\beta_k$.
Then, observe that \eqref{eq:DLMI_parameterization} can be expressed as
\begin{align}
(\star)^{\top}\underbrace{\left[\begin{array}{cc}
H_{k,k}^k  & \star\\
\frac{1}{2}\left(H_{k,k+1}^k+H_{k+1,k}^k\right) & H_{k+1,k+1}^k
\end{array}\right]}_{\coloneqq P}\underbrace{\left[\begin{array}{c}
\sigma_{1}^k(t)I\\
\sigma_{2}^k(t)I
\end{array}\right]}_{\coloneqq \Sigma}\succeq0.\label{eq:DLMI_sufficient_2}
\end{align}
The condition \eqref{eq:DLMI_sufficient_2} can be achieved by satisfying
$P\in\mtc$ where 
\begin{align*}
\mtc\coloneqq\left\{ P\bigg|\forall\sigma_1,\forall\sigma_2\in\mathbb{R}_{+},\Sigma=[\sigma_{1}I,\sigma_{2}I]^{\top},\Sigma^{\top}P\Sigma\succeq0\right\} .
\end{align*}
It is worth noting that if $n_{H}=1$, the set $\mtc$ is equal to
the set of \( 2 \times 2 \) copositive matrices \citep{boyd2004convex}. With this observation,
we are ready to derive sufficient LMI conditions ensuring $P\in\mtc$ by adapting techniques from the literature on matrix copositivity \citep{parrilo2000semidefinite, arceo2020copositive}. It is worth noting that we formulate these conditions specifically for our problem so that we can obtain a finite set of LMIs that guarantee the invariance of the funnel.

\begin{lemma}
\label{lem:first_copositive1}
Suppose that the following holds:
\begin{subequations}
\label{eq:LMI_copositive1}
\begin{align}
H_{k,k}^k & \succeq0,\quad H_{k+1,k+1}^k  \succeq0,\\
H_{k,k+1}^k&+H_{k+1,k}^k  \succeq0.
\end{align}
\end{subequations}
Then, \eqref{eq:DLMI_sufficient_2} holds for all $t$ in $(t_k,t_{k+1})$.
\end{lemma}
\begin{proof}
Each block of $P$ is PD because of the hypothesis in the Lemma \ref{lem:first_copositive1}. Since each $\sigma_{i}^k(t)$
for $i\in\{1,2\}$ is nonnegative for all $t$ in $(t_k,t_{k+1})$, we have \eqref{eq:DLMI_parameterization}.
\end{proof}
Similar results were discussed in \citep{oliveira2005stability,reynolds2021funnel}. Here, we can further derive a less conservative LMI condition implying $P\in\mtc$ compared to \eqref{eq:LMI_copositive1}.

\begin{lemma}
\label{lem:first_copositive2} Suppose that there exist $H_{i,j}^k$ and $X_{i,j}^k$ for all $i,j\in\{k,k+1\}$ such that
\begin{subequations}
\label{eq:LMI_copositive2}
\begin{align}
\left[\begin{array}{cc}
H_{k,k}^k-X_{k,k}^k  & \star\\
\frac{1}{2}\left(H_{k,k+1}^k+H_{k+1,k}^k\right)-X_{k+1,k}^k & H_{k+1,k+1}^k-X_{k+1,k+1}^k
\end{array}\right] & \succeq0,\\
X_{i,j}^k=(X_{j,i}^k)^{\top}\succeq0,\quad\forall i, j\in\{k,k+1\}.
\end{align}
\end{subequations}
Then, \eqref{eq:DLMI_sufficient_2} holds for all $t$ in $(t_k,t_{k+1})$.
\end{lemma}
\begin{proof}
Notice that the matrix $P$ defined in \eqref{eq:DLMI_sufficient_2}
can be written as
\begin{align*}
P & =P_{1}+P_{2},
\end{align*}
where $P_1$ and $P_2$ are given by
\begin{align*}
P_{1} & =\left[\begin{array}{cc}
H_{k,k}^k-X_{k,k}^k  & \star\\
\frac{1}{2}\left(H_{k,k+1}^k+H_{k+1,k}^k\right)-X_{k+1,k}^k & H_{k+1,k+1}^k-X_{k+1,k+1}^k
\end{array}\right],\\
P_{2} &= \left[\begin{array}{cc}
X_{k,k}^k  & \star\\
X_{k+1,k}^k  & X_{k+1,k+1}^k
\end{array}\right].
\end{align*}
We have $\Sigma^{\top}P_{1}\Sigma\succeq0$ because of the hypothesis in Lemma \ref{lem:first_copositive2}. Also, since every block $X_{i,j}^k$ for all $i,j\in\{k,k+1\}$ in $P_2$ is assumed to be PD, we have $\Sigma^{\top}P_{2}\Sigma\succeq0$. Hence, we can conclude $\Sigma^{\top}P \Sigma\succeq0$.
\end{proof}
By solving either \eqref{eq:LMI_copositive1}
or \eqref{eq:LMI_copositive2} for all $k$ in $\mathbb{Z}_{[0,N-1]}$, we can guarantee the satisfaction of the DLMI \eqref{eq:DLMI}
for all $t\in[t_{0},t_{f}]$ except at each temporal node points $t=t_k$ for all $k$ in $\mathbb{Z}_{[0,N]}$. 
\begin{rem}
    The matrix-valued function $Q(t)$ with the FOH \eqref{eq:FOH} is not differentiable at each temporal node points $t=t_k$ for all $k$ in $\mathbb{Z}_{[0,N]}$. This does not compromise our objective of continuous-time invariance since it suffices for the Lyapunov condition \eqref{eq:Lyapunov_condition} to hold for almost all $t$ in $[t_0,t_f]$ according to Lemma \ref{lem:invariance}.
\end{rem}

It is worth noting that \eqref{eq:LMI_copositive2} is less conservative than
\eqref{eq:LMI_copositive1}, as every solution of \eqref{eq:LMI_copositive1}
is a special case of \eqref{eq:LMI_copositive2}. This follows from the fact that  \eqref{eq:LMI_copositive2} with the following constraints
\begin{align*}
H_{i,i}^k & =X_{i,i}^k,\quad\forall i\in\{k,k+1\},\\
\frac{1}{2}\left(H_{k,k+1}^k+H_{k+1,k}^k\right) & =X_{k,k+1}^k,
\end{align*}
is equivalent to \eqref{eq:LMI_copositive1}.
However, solving \eqref{eq:LMI_copositive2} requires
to introduce additional variables $X\in\mathbb{R}^{n_{H}\times n_{H}}$
that are larger in dimension than the variables $Q$ and $Y$, making the solution of \eqref{eq:LMI_copositive2} more computationally expensive than that of \eqref{eq:LMI_copositive1}.

\subsection{Constraints}
The proposed funnel synthesis algorithm aims to satisfy not only the
invariance of the funnel but also the state and input
constraints. 
We consider linear state and input constraints written
as
\begin{subequations}
\label{eq:polytopoic_constraints}
\begin{align}
\mtp_{x} & =\{x\mid(a_{i}^{h})^{\top}x\leq b_{i}^{h},\quad i=1,\ldots,m_{x}\}, \\
\mtp_{u} & =\{x\mid(a_{j}^{g})^{\top}u\leq b_{j}^{g},\quad j=1,\ldots,m_{u}\}.\label{eq:input_constraint}
\end{align}
\end{subequations}
We aim to make the state and input funnels at each node point centered
around the nominal trajectory remain inside the feasible region. This
can be stated using the set inclusions as $\{\bar{x}_{k}\}\oplus\mathcal{E}_{Q}(t_{k})\subseteq\mtp_{x}$ and $\{\bar{u}_{k}\}\oplus\mte_{KQK^{\top}}(t_{k})\subseteq\mtp_{u}$
for all $k\in\mathbb{Z}_{[0,N]}$ where $\bar{x}_{k} = \bar{x}(t_k)$ and $\bar{u}_{k} = \bar{u}(t_k)$. These conditions can be equivalently
written by the following LMIs \citep{10167750}:
\begin{subequations}
\label{eq:funnel_constraints}
\begin{align}
0 & \preceq\left[\begin{array}{cc}
\left(b_{i}^{h}-(a_{i}^{h})^{\top}\bar{x}_{k}\right)^{2} & (a_{i}^{h})^{\top}Q_{k}\\
Q_{k}a_{i}^{h} & Q_{k}
\end{array}\right], i=\mathbb{Z}_{[1,m_x]},\\
0 & \preceq\left[\begin{array}{cc}
\left(b_{j}^{g}-(a_{j}^{g})^{\top}\bar{u}_{k}\right)^{2} & (a_{j}^{g})^{\top}Y_{k}\\
Y_{k}^{\top}a_{j}^{g} & Q_{k}
\end{array}\right], j=\mathbb{Z}_{[1,m_u]}.
\end{align}
\end{subequations}

\subsection{Problem formulation}

The proposed funnel synthesis algorithm solves the following SDP problem:
\begin{subequations}
\label{eq:SDP}
\begin{align}
\underset{\resizebox{0.33\hsize}{!}{$Q_{k},Y_{k},(\lambda_\beta)_{k},(\lambda_\gamma)_{k},X_{k,k}^{k},X_{k,k+1}^{k},X_{k+1,k+1}^{k}$}}{\operatorname{minimize}}~& J(Q_{0},\ldots,Q_{N})\\
 \operatorname{subject~to}~~~~& \eqref{eq:LMI_copositive1}\text{ or }\eqref{eq:LMI_copositive2},\forall k\in\mathbb{Z}_{[0,N-1]},\label{eq:LMI_FOH_in_SDP}\\
 & \text{\eqref{eq:funnel_constraints}},\forall k\in\mathbb{Z}_{[0,N]},\\
 & 0\prec Q_{k}\preceq Q_{max},\forall k\in\mathbb{Z}_{[0,N]}.
\end{align}
\end{subequations}
The cost function $J$ is assumed to be convex in $Q_{k}$. 
The matrix $Q_{max}\in\mathbb{S}_{++}^{n_{x}}$ is introduced to prohibit
the funnel being arbitrarily large. By imposing $0\prec Q_{k}$ for
each $k\in\mathbb{Z}_{[0,N]}$, the FOH interpolation \eqref{eq:FOH}
ensures $0\prec Q(t)$ for all $t\in[t_{0},t_{f}]$.

\begin{thm}
\label{thm:final}
Let $Q(t)$ and $Y(t)$ be obtained from the solution of \eqref{eq:SDP} with the FOH interpolation \eqref{eq:FOH}. Define the feedback gain $K(t) = Y(t) Q(t)^{-1}$. Then, the funnel $\mathcal{F}(t)$ defined in \eqref{eq:deffunnel} is invariant for all $t \in [t_0, t_f]$.
\end{thm}

\begin{proof}
By Lemma \ref{lem:first_copositive1} and \ref{lem:first_copositive2}, satisfying \eqref{eq:LMI_FOH_in_SDP}
implies that the DLMI \eqref{eq:DLMI} holds for all $t\in[t_{0},t_{f}]$
except at each temporal node point $t=t_{k}$ for all $k\in\mathbb{Z}_{[0,N]}$.
It follows from Lemma \ref{lem:DLMI} that the Lyapunov condition
\eqref{eq:Lyapunov_condition} holds almost everywhere. This implies
the invariance of the state funnel $\mte_{Q}$ followed from Lemma \ref{lem:invariance}. Last, it follows from Lemma \ref{lem:nonlinear_invariance} that the funnel $\mtf(t)$ is invariant.
\end{proof}
\begin{rem}
\label{rem:CTCS} The feasibility condition in \eqref{eq:funnel_constraints}
enforces the constraints only at the temporal node points, which could
result in violations between these node points. A possible remedy
is to increase the number of node points. It is important to note
that while the feasibility is only enforced at each node point, the
invariance holds across the entire horizon $[t_{0},t_{f}]$ as established
in Theorem \ref{thm:final}.
\end{rem}
\begin{rem}
\label{rem:infeasibility} 
The feasibility of the problem \eqref{eq:SDP} depends on many factors such as the magnitude of $\gamma$ and $\beta$, the density of the node points, the configuration and size of obstacles, and the bounds on control inputs. Due to this complexity, it is generally difficult to guarantee feasibility a priori. However, a key advantage of the proposed approach is that it is formulated as a convex optimization problem, certifying the problem's infeasibility. In particular, if the problem is infeasible, the solver can provide a formal certificate of infeasibility.
\end{rem}

\subsection{Estimation of the constants via sampling}

To solve the problem \eqref{eq:SDP}, it is necessary to obtain the
local Lipschitz constant $(\gamma_i)_{k}$ and the bound $\beta_{k}$
described in \eqref{eq:bound_on_Delta} for each $k\in\mathbb{Z}_{[0,N-1]}$.
If the system is globally Lipschitz, one can directly use the global Lipschitz constant for all $k$ \citep{yu2013tube}. One effective way for estimating local $\gamma_{k}$ is to use sampling
approaches \citep{reynolds2021funnel}, collecting samples of the
triple $(\eta_{s},\xi_{s},w_{s})$ around the nominal trajectory,
with $s$ denoting each sample. Specifically, this involves sampling
the state deviation $\eta_{s}$ within the maximum state funnel $Q_{max}$,
the input deviation $\xi_{s}$ within the constraint \eqref{eq:input_constraint},
and the disturbance within its bound \eqref{eq:bound_on_2} for each
subinterval $[t_{k},t_{k+1}]$. 
To obtain a less conservative constant, one might use the iterative procedure provided in \citep{reynolds2021funnel} where the estimation of $\gamma_{k}$ and funnel computation are alternately repeated until the convergence. This paper does not consider this iterative
process for simplicity as our proposed
method can be incorporated with any way to estimate the Lipschitz constant. On the other hand, the value of $\Delta_e(t)$
defined in \eqref{eq:LTV_1} is known and can thus be evaluated at
each $t$. Hence, for each subinterval $[t_{k},t_{k+1}]$, we pick
$N_{s}$ sample temporal points $t_{s}$ uniformly, and then $(\beta_e)_{k}$
can be determined through 
\[
(\beta_e)_{k}=\max_{s=1,\ldots,N_{s}}\norm{\Delta_e(t_{s})}_{2}.
\]

\subsection{Algorithm summary}
The overall procedure for synthesizing the funnel and controller is summarized below.

\begin{algorithm}[H]
\caption{Funnel Synthesis}
\begin{algorithmic}[1]
\State \textbf{Input:} Nominal trajectory, system dynamics, state/input constraints
\State Estimate $\gamma$ and $\beta$ as described in Section~3.4
\For{each $\lambda_w$ in the search grid}
    \State Solve the SDP problem in (\ref{eq:SDP})
    \State Store the cost and feasibility status
\EndFor
\State Select $\lambda_w$ with the lowest feasible cost
\State Return the corresponding funnel and controller
\end{algorithmic}
\end{algorithm}

\section{Numerical simulation}

In this section, the proposed method is validated through two robotic motion planning and control applications: one for a unicycle model and the other for six-degree-of-freedom (6-DoF) powered descent guidance. For both examples, the Mosek solver is used to solve the SDP \eqref{eq:SDP}. 
The code
for the simulation can be found in the following repository.
\begin{center}
\vspace{-0.2cm}
{\footnotesize\href{https://github.com/taewankim1/funnel-copositivity}{\texttt{https://github.com/taewankim1/funnel-copositivity}}}
\vspace{-0.2cm}
\end{center}

\subsection{Unicycle model}

\begin{figure}
\begin{center}
\begin{subfigure}
 \centering
 \includegraphics[width=7.0cm,trim={0cm 0.0cm 0cm 0cm},clip]{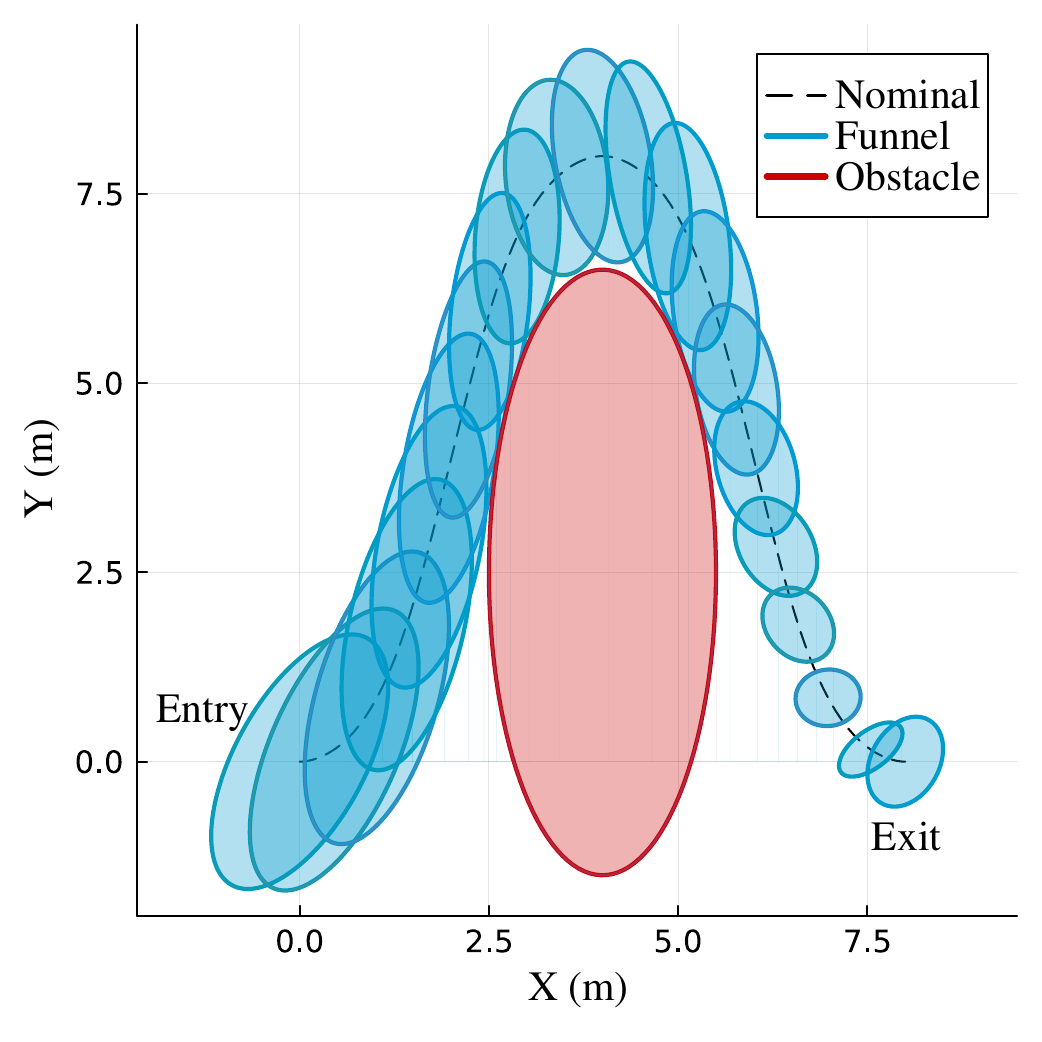}
\end{subfigure}
\begin{subfigure}
 \centering
 \includegraphics[width=8.0cm,clip]{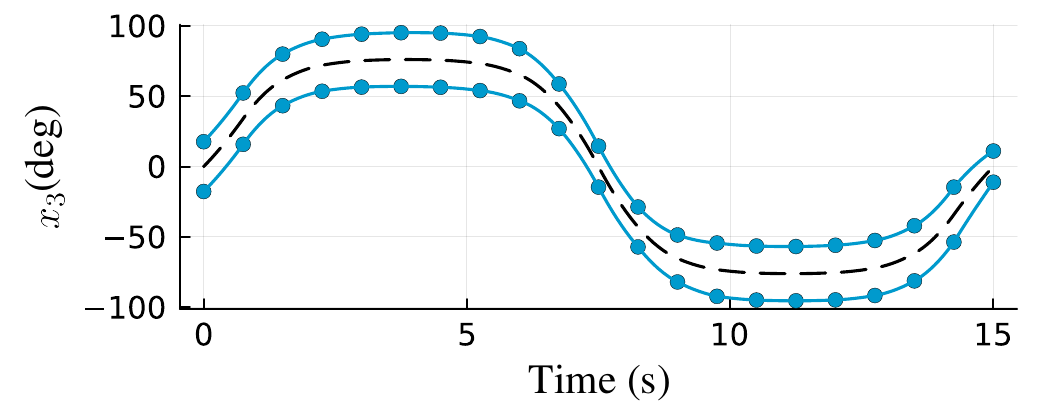}
\end{subfigure}
\caption{(Top) The synthesized state funnel projected on $x$ ($x_1$) and $y$ ($x_2$) position coordinates. (Bottom) Time history of the state funnel projected on yaw angle ($x_3$) coordinate.}
\label{fig:statefunnel}
\end{center}
\end{figure}

We consider the unicycle model given in \eqref{eq:unicycle} with
$c_{1}=0.03$ and $c_{2}=0.05$. We divide a $15$ second time horizon
into $N=20$ subintervals, starting at $t_{0}=0$ and ending at $t_{f}=15$.
The cost function $J$ in \eqref{eq:SDP} is set as $-\trace{Q_{0}}+\trace{Q_{N}}$
to maximize the funnel entry and minimize the funnel exit. The nominal
trajectory illustrated in Fig. \ref{fig:statefunnel} starts at $(0, 0)$, passes through $(4, 8)$, and ends at $(8, 0)$. A state constraint for avoiding an obstacle is taken into account, resulting in a nonconvex state constraint as
shown in Fig. \ref{fig:statefunnel}. This nonconvex constraint
is linearized around the nominal trajectory to generate the polytopic
constraints in the form of \eqref{eq:polytopoic_constraints}. The
constraints for inputs are defined as follows: $0\leq u_{1}\leq2$
(m/s) and $\abs{u_{2}}\leq1.5$ (rad/s). The parameter $Q_{max}$ for the maximum funnel
size is set to $\text{diag}(2^{2},2^{2},(20\pi/180)^{2})$. To estimate
$(\gamma_i){k}$ and $\beta_{k}$, we sample 100 triples of $(\eta_{s},\xi_{s},w_{s})$
and 20 samples of $\Delta(t_{s})$ for each subinterval $[t_{k},t_{k+1}]$.

\begin{figure}
\centering
\includegraphics[width=10cm,clip]{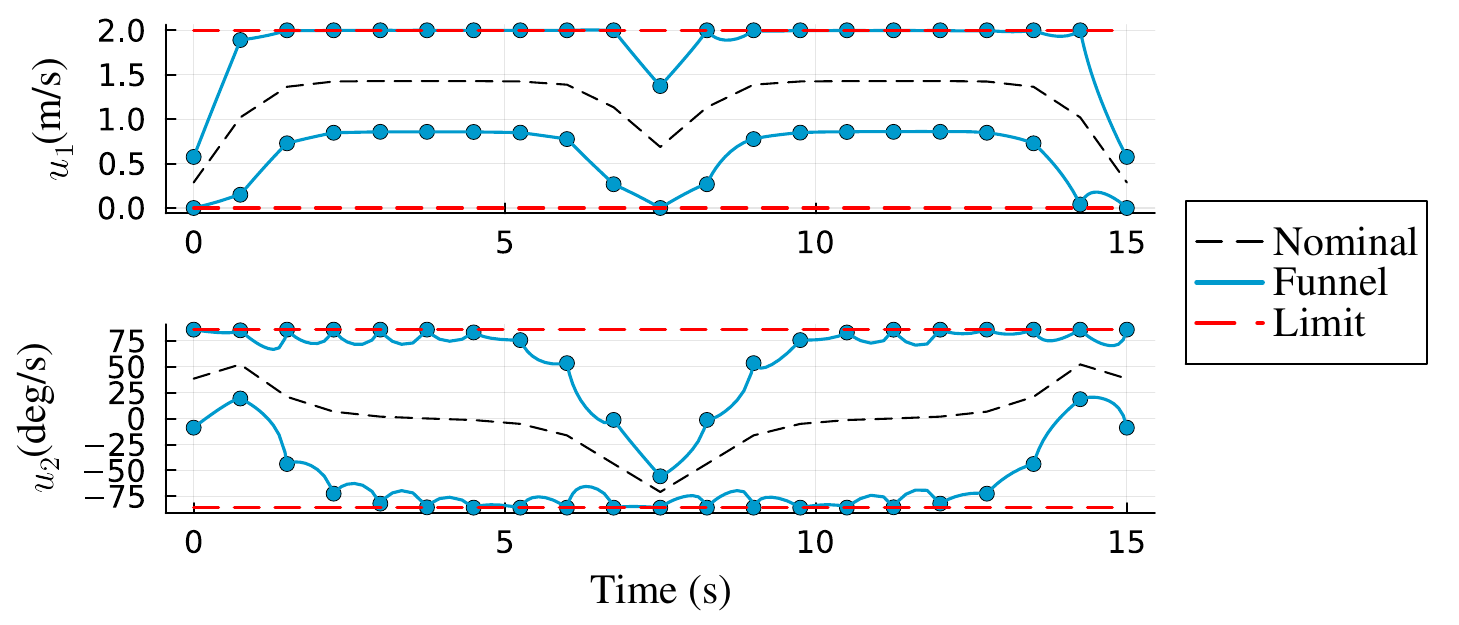}
\caption{Time history of the synthesized input funnel projected on velocity command ($u_1$) and angular velocity command ($u_2$) coordinates, shown in the top and bottom figures, respectively.}
\label{fig:inputfunnel}
\end{figure}

\begin{figure}
\centering
\includegraphics[width=8cm,clip]{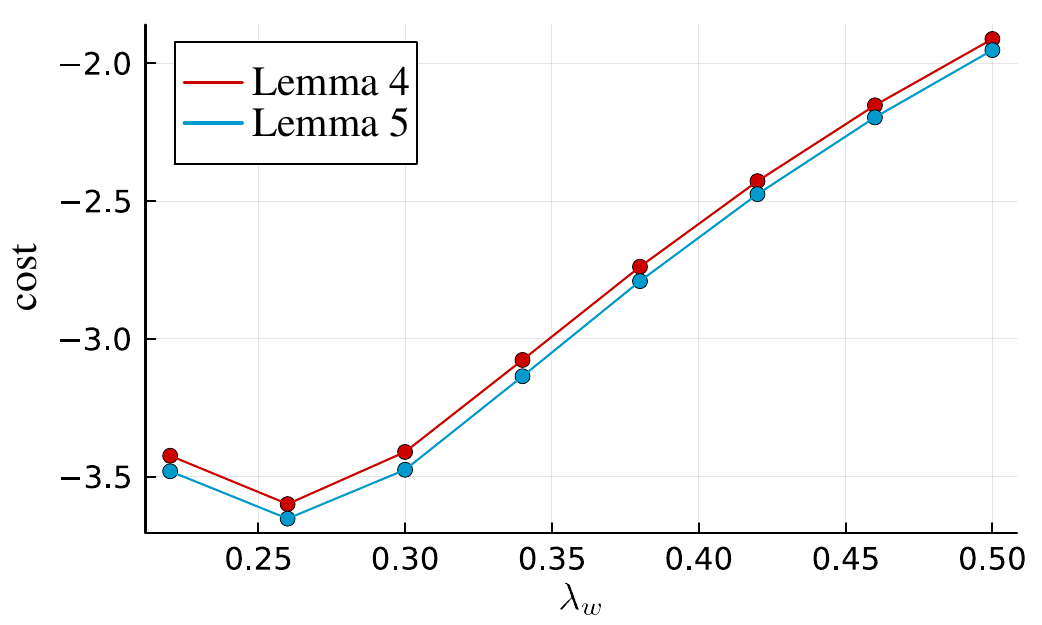}
\caption{The cost results of \eqref{eq:SDP} with different values of $\lambda_w$ are presented for both cases: using \eqref{eq:LMI_copositive1} in Lemma~\ref{lem:first_copositive1} and \eqref{eq:LMI_copositive2} in Lemma~\ref{lem:first_copositive2}.}
\label{fig:lambda_w_comparison}
\end{figure}

\begin{figure}
\begin{center}
\begin{subfigure}
 \centering
 \includegraphics[width=7.0cm,trim={0cm 0.0cm 0cm 0cm},clip]{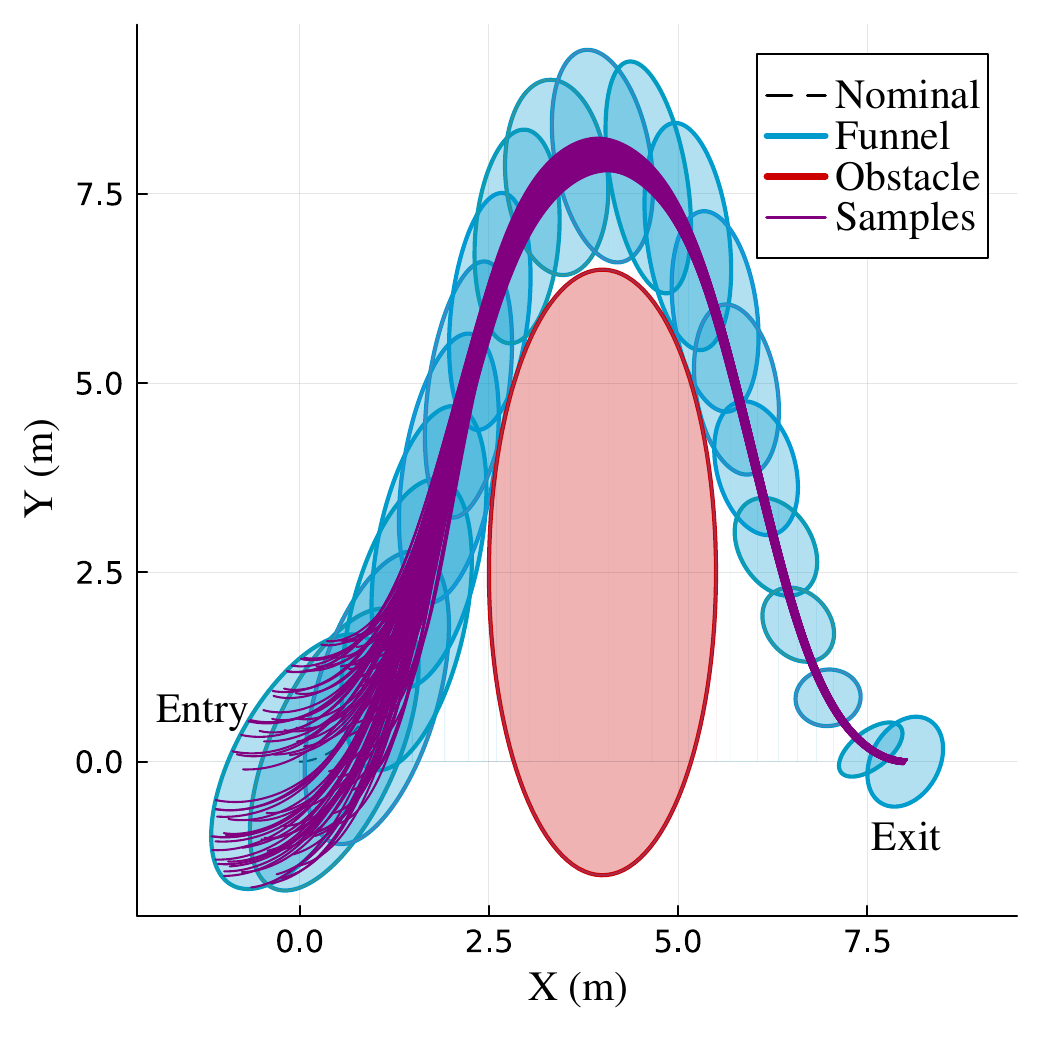}
\end{subfigure}
\begin{subfigure}
 \centering
 \includegraphics[width=8.0cm,clip]{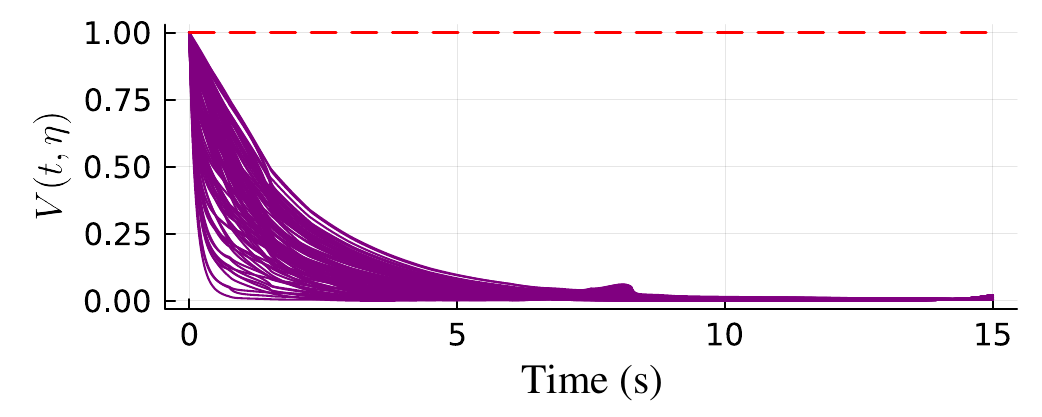}
\end{subfigure}
\caption{(Top) The results of trajectories propagated from randomly selected samples within the funnel entry. (Bottom) Time history of Lyapunov function $V$, as defined in \eqref{eq:Lyapunov_function}, for each trajectory sample.}
\label{fig:samples}
\end{center}
\end{figure}

\begin{figure}
\begin{center}
\begin{subfigure}
 \centering
 \includegraphics[width=8.0cm,clip]{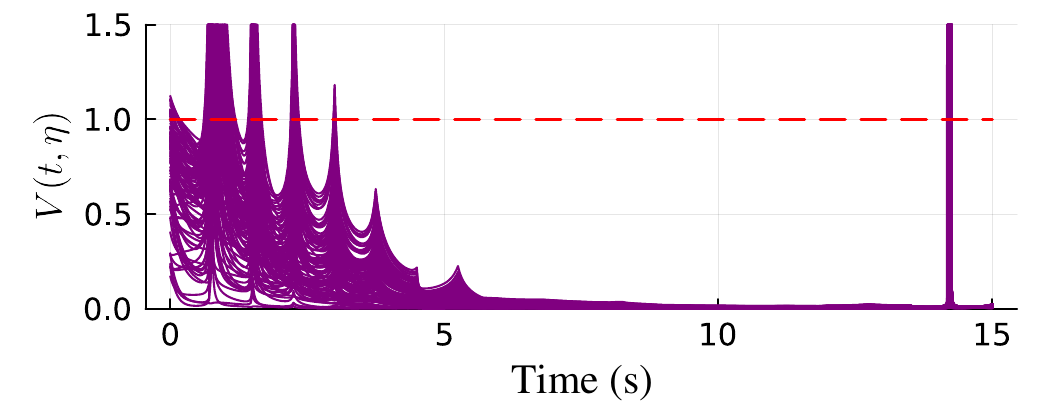}
\end{subfigure}
\begin{subfigure}
 \centering
 \includegraphics[width=8.0cm,clip]{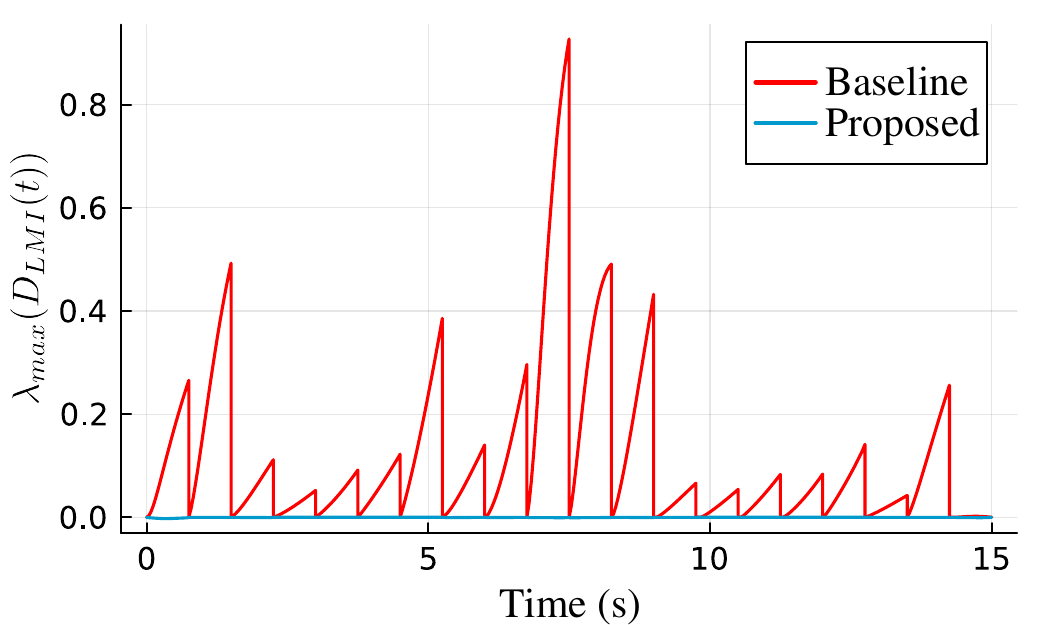}
\end{subfigure}
\caption{(Top) Evolution of the Lyapunov function for the funnel computed using the baseline approach. (Bottom) Time history of the maximum eigenvalue of $D_{LMI}$, as defined in \eqref{eq:DLMI}, for the funnel obtained by both the baseline and proposed methods.}
\label{fig:unicycle_baseline_samples}
\end{center}
\end{figure}

To determine a value for $\lambda_w$, we solved the SDP problem \eqref{eq:SDP} across a range of candidate $\lambda_w$ values. The resulting cost values are summarized in Fig. \ref{fig:lambda_w_comparison}, which shows that the minimum cost is attained when $\lambda_w = 0.26$. The results of the state funnel and input funnel are provided in Fig.
\ref{fig:statefunnel} and Fig. \ref{fig:inputfunnel}, respectively. It is clear that state and input constraints are satisfied at each
node point. We conduct a comparison between two different conditions,
\eqref{eq:LMI_copositive1} and \eqref{eq:LMI_copositive2}, in terms of the resulting
cost and computation time for synthesizing the funnel. The funnel
computed using \eqref{eq:LMI_copositive2} has smaller cost, at -3.652, compared
to -3.600 for the funnel using \eqref{eq:LMI_copositive1}. The
computational time for MOSEK to solve the SDP \eqref{eq:SDP} with \eqref{eq:LMI_copositive1}
is longer, at 11.63s, compared to 1.57s with \eqref{eq:LMI_copositive2},
which is aligned with our expectation. Furthermore, Fig. \ref{fig:lambda_w_comparison} shows that for every case, the cost associated with \eqref{eq:LMI_copositive2} is consistently lower than that obtained by \eqref{eq:LMI_copositive1}.

To test the invariance condition of the funnel, we take a
total of 500 samples from the surface of the funnel entry $Q_{0}$.
Then, these samples are propagated through the control law \eqref{eq:control_law}.
During the propagation, we randomly choose the disturbance $w$ such
that $\norm w_{2}=1$. In the bottom figure of Fig. \ref{fig:samples},
the value of the Lyapunov function $V$ for each sample over the time
horizon is illustrated. We can see that the values maintain below
one, which means all trajectories remain in the funnel by the invariance
property.

We compare the proposed method with a baseline approach that directly discretizes the DLMI \eqref{eq:DLMI} using the forward Euler method. This type of discretization has been commonly used in related funnel computation literature \cite[Section 4.2]{majumdar2017funnel}. To validate the baseline approach, we follow the same procedure as with the proposed method: 500 samples are taken from the surface of the funnel entry computed using the baseline, and each sample is propagated under the control law. The evolution of the Lyapunov function for each sample is shown in the top of Fig.~\ref{fig:unicycle_baseline_samples}. Some trajectories exhibit Lyapunov function values greater than 1, indicating that the funnel computed by the baseline method violates the invariance condition. To explicitly assess satisfaction of $\eqref{eq:DLMI}$ over time, we plot the time history of the maximum eigenvalue of $D_{LMI}$, which is defined in \eqref{eq:DLMI}, using the decision variables from both the baseline and proposed methods in Fig. \ref{fig:unicycle_baseline_samples}. The results clearly show that the baseline approach violates the DLMI condition between the node points, whereas the proposed method satisfies it throughout the entire time horizon.

\subsection{6-DoF Powered Descent Guidance}

We consider a 6-DoF rigid-body rocket model in an East-North-Up inertia coordinate described by
\begin{align*}
    \left[\begin{array}{c}
    \dot{p}\\
    \dot{v}\\
    \dot{\Phi} \\
    \dot{\Omega}
    \end{array}\right] & =\left[\begin{array}{c}
    v\\
    \frac{1}{m} C_{I/B} (\Phi) F + g\\
    R(\Phi) \Omega \\
    J^{-1} (T + r \times F - \Omega \times J\Omega)
    \end{array}\right],
\end{align*}
where the state vector $x = (p,v,\Phi,\Omega)$ consists of the position $p\in \mathbb{R}^3$, velocity $v\in \mathbb{R}^3$, Euler angles $\Phi\in \mathbb{R}^3$, and angular velocity $\Omega\in \mathbb{R}^3$, and the control input $u=(F,T)$ includes the thrust force $F\in \mathbb{R}^3$ and torque $T\in \mathbb{R}^3$. The constant $m \in \mathbb{R}$ is the vehicle mass, $J \in \mathbb{R}^{3\times3}$ is the inertia matrix, and $g\in\mathbb{R}^3$ is the gravitational constant. The vector $r\in\mathbb{R}^3$ is the position vector from the center of the mass to the location of the thrust. The matrix $C_{I/B}(\Phi)$ is the  rotation matrix from the body frame to the inertia frame and the matrix $R(\Phi)$ represents the transformation between Euler angle rates and body-frame angular velocity. More details in the derivation of the model could be found in \cite{doi:10.2514/1.G004549}.

The nominal trajectory has a plane maneuver in the plane defined by $r_x=r_y$, with $r_x$ and $r_y$ denoting the position components along the East and North axes, respectively. We use nondimensionalized units, with $U_T$, $U_L$, and $U_M$ denoting the units of time, length, and mass units, respectively. The mass $m$ is set to $2$ $(U_M)$ and the gravity is set to $g=[0, 0, -1.625]^\top$ $(U_M/U_T^2)$. We divide an entire horizon into $N = 5$ subintervals, starting at $t_0 = 0$ $(U_T)$ and ending at $t_f = 7.457$ $(U_T)$. The cost function is set to the same as in the unicycle example: $-\trace{Q_{0}}+\trace{Q_{N}}$. The constraint sets $\mathcal{X}$ and $\mathcal{U}$ are defined as
\begin{align*}
    \mathcal{X} &= \{x \mid x_{lb} \leq x \leq x_{ub} \},\\
    \mathcal{U} &= \{u \mid u_{lb} \leq u \leq u_{ub} \}, \\
    x_{ub} &= (8,8,8,2,2,2,30\degree,30\degree,30\degree,45\degree,45\degree,45\degree), \\
    x_{lb} &= -x_{lb}, \\
    u_{ub} &= (1.5,1.5,5.0,0.1,0.1,0.1), \\
    u_{lb} &= -(1.5,1.5,0.0,0.1,0.1,0.1).
\end{align*}
Since the bounded disturbance $w$ is not considered in this example, it is not necessary to perform the line search for $\lambda_w$; thus $\lambda_w$ is set to zero and the matrix $F(t)$ is set to the zero matrix for all time.

\begin{figure}
\centering
\includegraphics[width=8cm,clip]{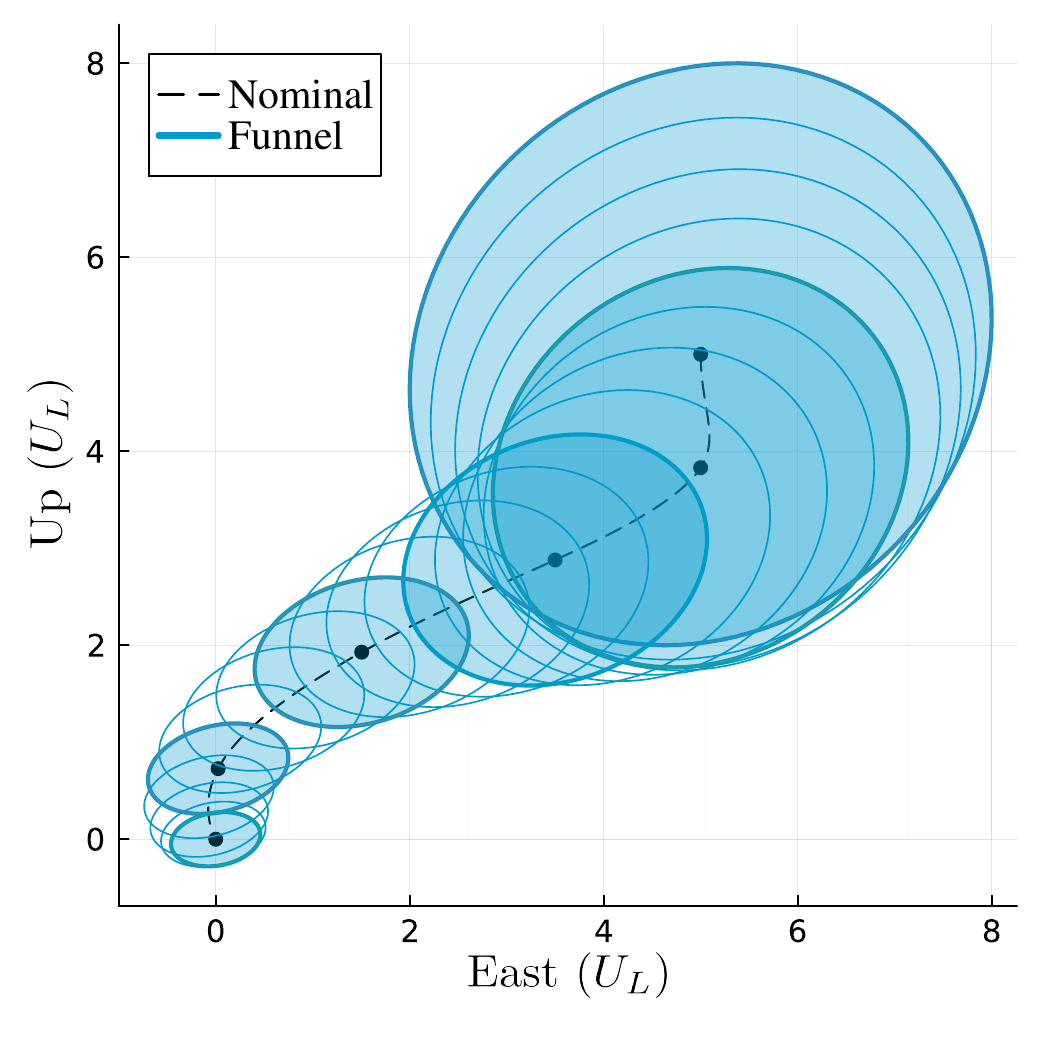}
\caption{The synthesized state funnel projected onto the position coordinates. Filled ellipsoids represent the funnels at discrete node points, while the intermediate funnels between nodes are shown as unfilled ellipsoids.}
\label{fig:PDG_state_funnel}
\end{figure}

\begin{figure}
\begin{center}
\begin{subfigure}
 \centering
 \includegraphics[width=8.0cm,clip]{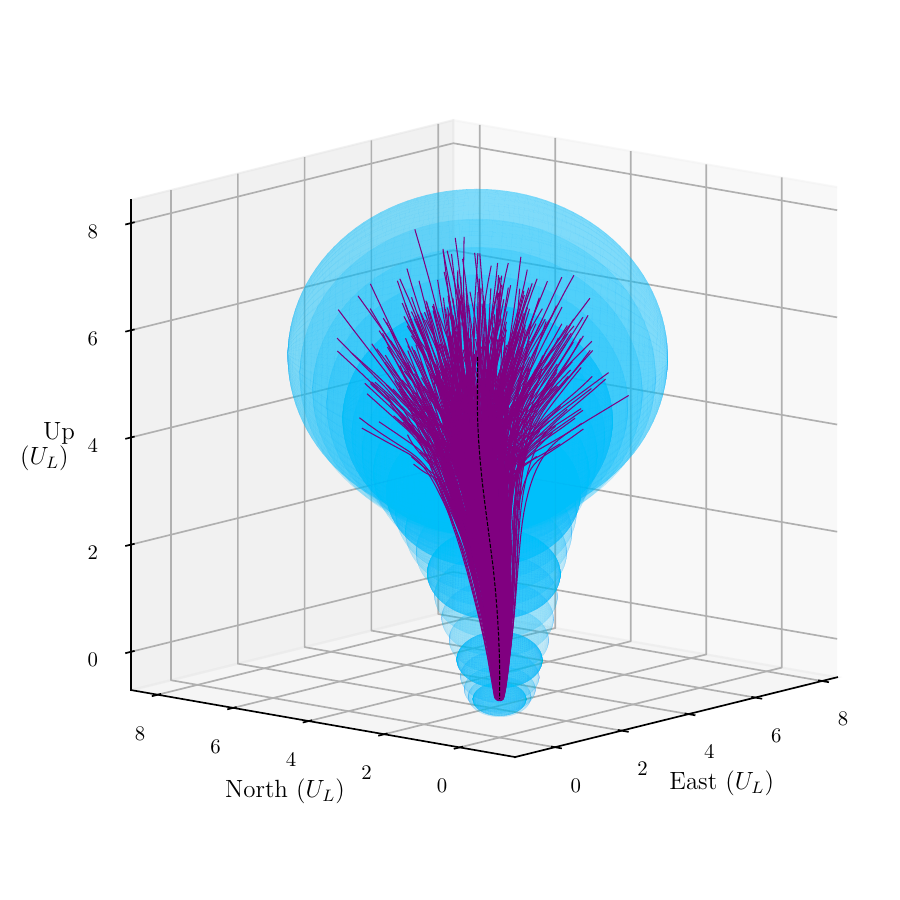}
\end{subfigure}
\begin{subfigure}
 \centering
 \includegraphics[width=8.0cm,clip]{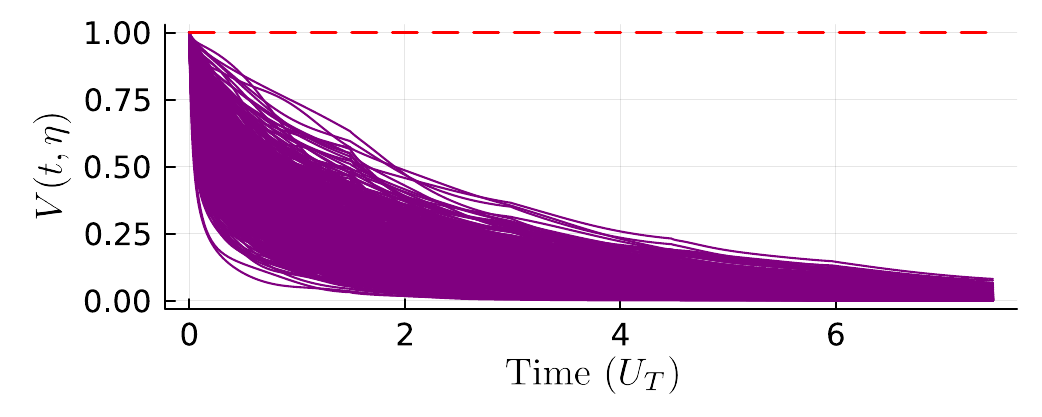}
\end{subfigure}
\caption{(Top) The computed funnel projected onto the three-dimensional position space. The funnel is shown in blue, and the propagated sample trajectories are shown in purple. (Bottom) Evolution of the Lyapunov function of the sample trajectories.}
\label{fig:PDG_samples}
\end{center}
\end{figure}

The computed funnels projected onto the two dimensional East-North plane and the three-dimensional East-North-Up space are illustrated in Fig. \ref{fig:PDG_state_funnel} and in Fig. \ref{fig:PDG_samples}, respectively. The funnel computed using \eqref{eq:LMI_copositive2} has a cost of -31.982, which is lower than the cost of the funnel computed using \eqref{eq:LMI_copositive1}, -31.890. The MOSEK solver takes 28.10 seconds to solve the SDP with \eqref{eq:LMI_copositive2} and 17.41 seconds with \eqref{eq:LMI_copositive1}.

\begin{figure}
\begin{center}
\begin{subfigure}
 \centering
 \includegraphics[width=8.0cm,clip]{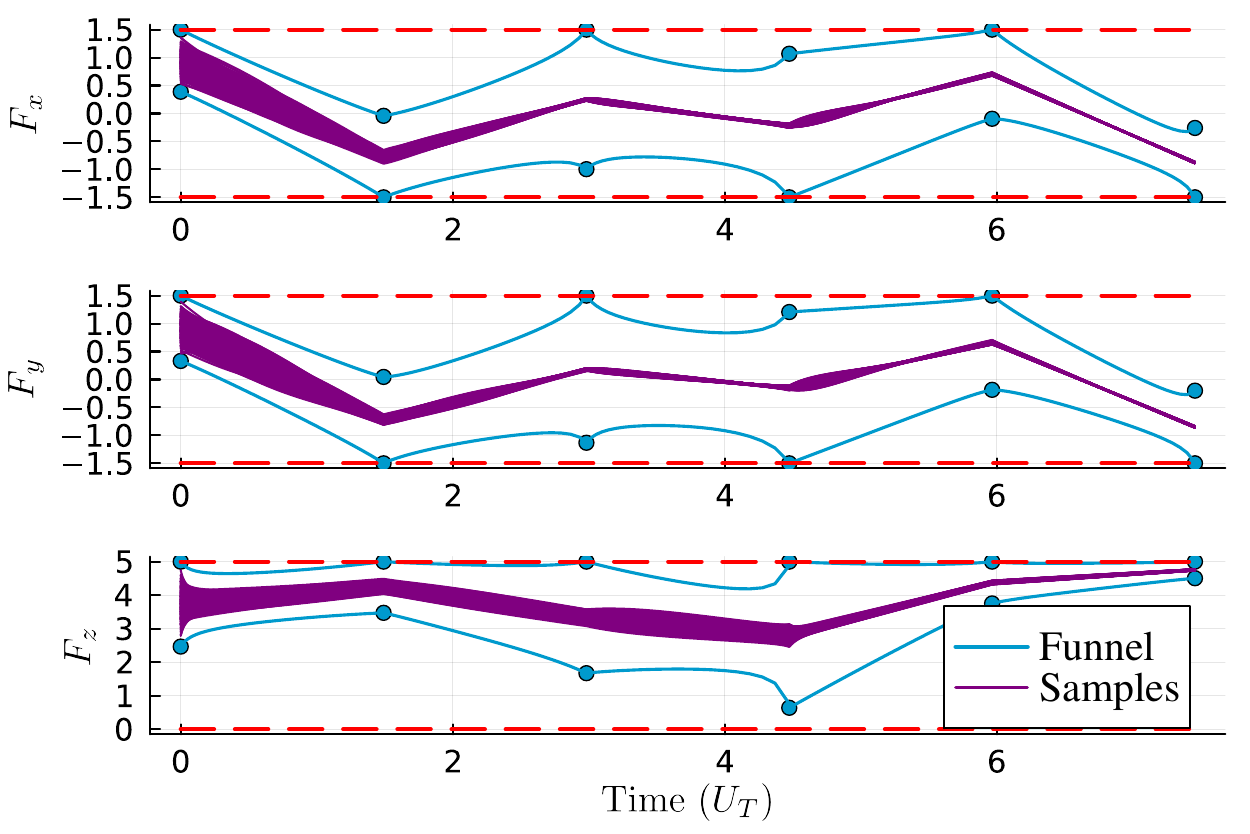}
\end{subfigure}
\begin{subfigure}
 \centering
 \includegraphics[width=8.0cm,clip]{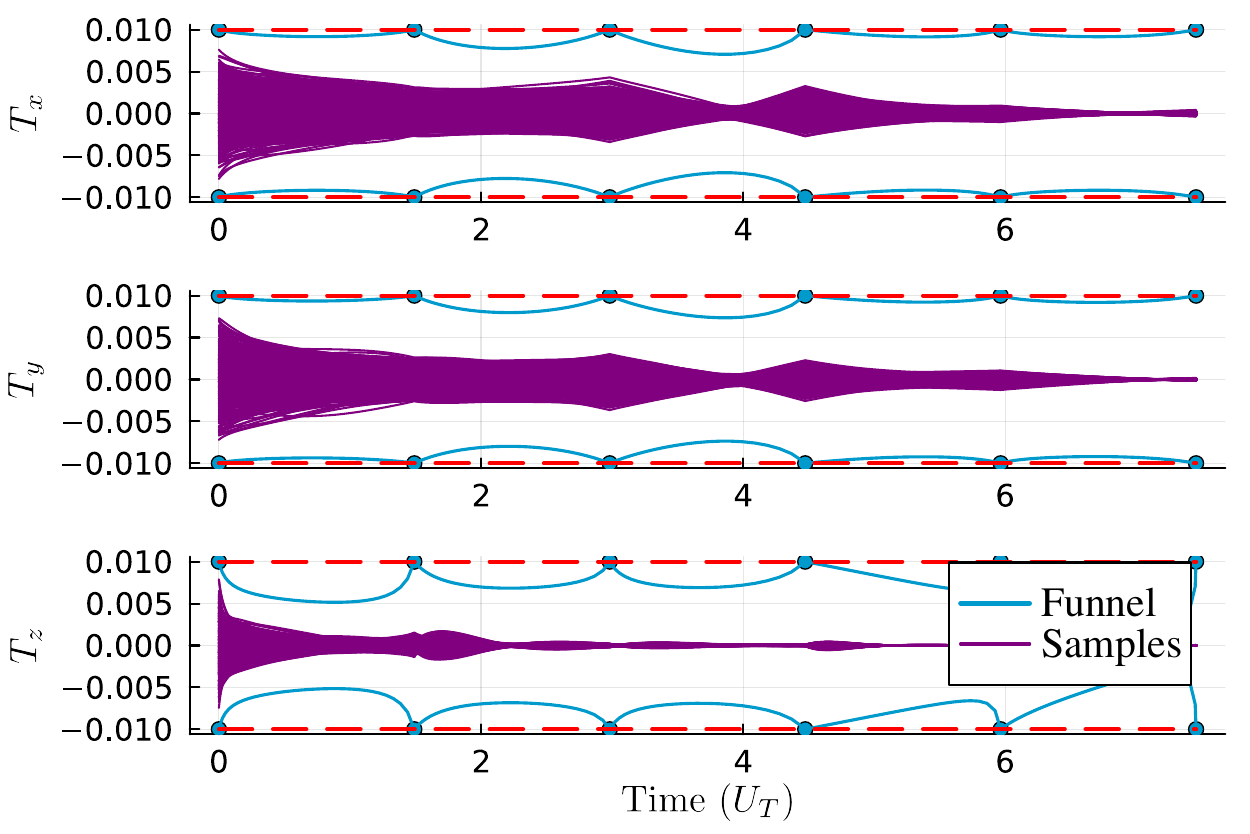}
\end{subfigure}
\caption{The synthesized input funnel onto each input dimension.}
\label{fig:PDG_input_funnels}
\end{center}
\end{figure}

\begin{figure}
\centering
\includegraphics[width=8cm,clip]{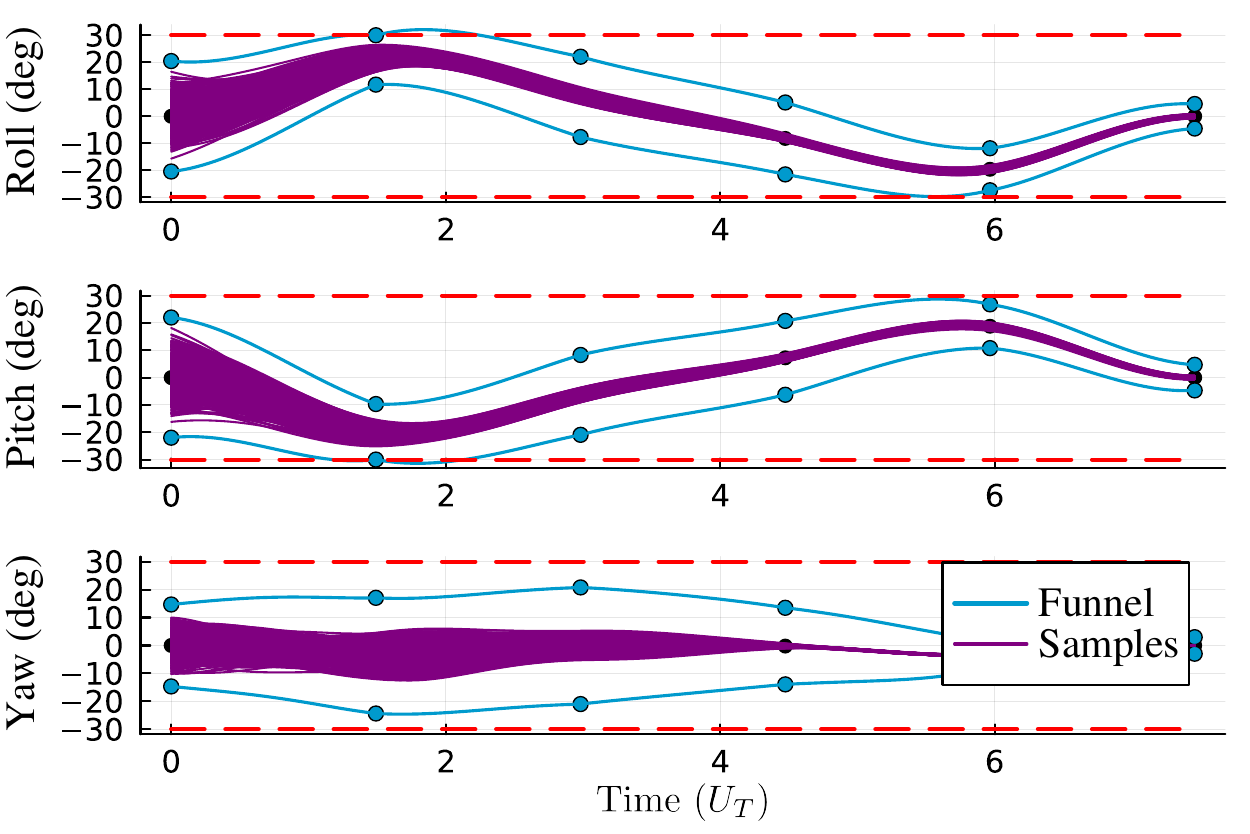}
\caption{The synthesized state funnel projected onto each Euler angle dimension.}
\label{fig:PDG_euler_angles}
\end{figure}

Similar to the unicycle example, 500 samples are taken from the funnel entry, and each sample is propagated under the synthesized control law. Each sample's trajectory in the position coordinates is illustrated in Fig. \ref{fig:PDG_samples}. From the results, we observe that although the nominal trajectory has an in-plane maneuver, the resulting funnel and sampled trajectories exhibit out-of-plane maneuver. The funnel invariance is validated by examining the evolution of the Lyapunov function, as shown in Fig. \ref{fig:PDG_euler_angles}, confirming that all trajectories remain inside the funnel. The resulting input funnel, projected onto each input dimension, is given in Fig. \ref{fig:PDG_input_funnels}, demonstrating that the input constraints are satisfied. The state funnel projected onto the Euler angle dimensions is given in Fig. \ref{fig:PDG_euler_angles}. While the funnel satisfies the state constraints at each node point, slight constraint violations occur between nodes.

\section{Conclusion}

This paper presents a funnel synthesis algorithm for nonlinear
systems under bounded disturbances. The proposed method can satisfy
the invariance condition exactly over the finite time horizon. This
is achieved by approximating the incremental system as an uncertain
LPV system and deriving a finite number of LMIs that imply the DLMI
resulting from the invariance condition. 
Simulation results show that the synthesized funnel satisfies both invariant and feasibility conditions. Future work will focus on extending the proposed method to continuous-time constraint satisfaction, not only for the invariance condition but also for state and input constraints. Improving computational efficiency for large-scale SDP problems is also a valuable direction.

\section*{Acknowledgement}
This work was supported in part by the Air Force Office of
Scientific Research under Grant FA9550-20-1-0053 and in part by the Office of Naval Research under Grant ONR N000142512231.

\bibliographystyle{elsarticle-num} 
\bibliography{root}

\begin{thebibliography}{10}
\expandafter\ifx\csname url\endcsname\relax
  \def\url#1{\texttt{#1}}\fi
\expandafter\ifx\csname urlprefix\endcsname\relax\def\urlprefix{URL }\fi
\expandafter\ifx\csname href\endcsname\relax
  \def\href#1#2{#2} \def\path#1{#1}\fi

\bibitem{mayne2005robust}
D.~Q. Mayne, M.~M. Seron, S.~Rakovi{\'c}, Robust model predictive control of
  constrained linear systems with bounded disturbances, Automatica 41~(2)
  (2005) 219--224.

\bibitem{rakovic2016elastic}
S.~V. Rakovi{\'c}, W.~S. Levine, B.~A{\c{c}}ikmese, Elastic tube model
  predictive control, in: 2016 American Control Conference (ACC), IEEE, 2016,
  pp. 3594--3599.

\bibitem{majumdar2017funnel}
A.~Majumdar, R.~Tedrake, Funnel libraries for real-time robust feedback motion
  planning, The International Journal of Robotics Research 36~(8) (2017)
  947--982.

\bibitem{reynolds2021funnel}
T.~Reynolds, D.~Malyuta, M.~Mesbahi, B.~A{\c{c}}{\i}kme{\c{s}}e, J.~M. Carson,
  Funnel synthesis for the 6-{DoF} powered descent guidance problem, in: AIAA
  SciTech 2021 Forum, 2021, p. 0504.

\bibitem{seo2021fast}
H.~Seo, C.~Y. Son, H.~J. Kim, Fast funnel computation using multivariate
  {B}ernstein polynomial, IEEE Robotics and Automation Letters 6~(2) (2021)
  1351--1358.

\bibitem{10167750}
T.~Kim, P.~Elango, T.~P. Reynolds, B.~A{\c{c}}{\i}kme{\c{s}}e, M.~Mesbahi,
  Optimization-based constrained funnel synthesis for systems with {L}ipschitz
  nonlinearities via numerical optimal control, IEEE Control Systems Letters 7
  (2023) 2875--2880.
\newblock \href {https://doi.org/10.1109/LCSYS.2023.3290229}
  {\path{doi:10.1109/LCSYS.2023.3290229}}.

\bibitem{tobenkin2011invariant}
M.~M. Tobenkin, I.~R. Manchester, R.~Tedrake, Invariant funnels around
  trajectories using sum-of-squares programming, IFAC Proceedings Volumes
  44~(1) (2011) 9218--9223.

\bibitem{jang2021fast}
I.~Jang, H.~Seo, H.~J. Kim, Fast computation of tight funnels for piecewise
  polynomial systems, IEEE Control Systems Letters (2021).

\bibitem{buch2021finite}
J.~Buch, P.~Seiler, Finite horizon robust synthesis using integral quadratic
  constraints, International Journal of Robust and Nonlinear Control 31~(8)
  (2021) 3011--3035.

\bibitem{seiler2019finite}
P.~Seiler, R.~M. Moore, C.~Meissen, M.~Arcak, A.~Packard, Finite horizon
  robustness analysis of {LTV} systems using integral quadratic constraints,
  Automatica 100 (2019) 135--143.

\bibitem{seiler2024trajectory}
P.~Seiler, R.~Venkataraman, Trajectory-based robustness analysis for nonlinear
  systems, International Journal of Robust and Nonlinear Control 34~(2) (2024)
  910--926.

\bibitem{pfifer2015robustness}
H.~Pfifer, P.~Seiler, Robustness analysis of linear parameter varying systems
  using integral quadratic constraints, International Journal of Robust and
  Nonlinear Control 25~(15) (2015) 2843--2864.

\bibitem{parrilo2000semidefinite}
P.~A. Parrilo, Semidefinite programming based tests for matrix copositivity,
  in: Proceedings of the 39th IEEE Conference on Decision and Control (Cat. No.
  00CH37187), Vol.~5, IEEE, 2000, pp. 4624--4629.

\bibitem{arceo2020copositive}
J.~C. Arceo, J.~Lauber, Copositive conditions for {LMI}-based controller and
  observer design, IFAC-PapersOnLine 53~(2) (2020) 7959--7964.

\bibitem{d2013incremental}
L.~D'Alto, M.~Corless, Incremental quadratic stability, Numerical Algebra,
  Control and Optimization 3~(1) (2013) 175--201.

\bibitem{xu2020observer}
X.~Xu, B.~A{\c{c}}{\i}kme{\c{s}}e, M.~J. Corless, Observer-based controllers
  for incrementally quadratic nonlinear systems with disturbances, IEEE
  Transactions on Automatic Control 66~(3) (2020) 1129--1143.

\bibitem{accikmecse2011robust}
B.~A{\c{c}}{\i}kme{\c{s}}e, J.~M. Carson~III, D.~S. Bayard, A robust model
  predictive control algorithm for incrementally conic uncertain/nonlinear
  systems, International Journal of Robust and Nonlinear Control 21~(5) (2011)
  563--590.

\bibitem{malikov2020numerical}
A.~I. Malikov, D.~I. Dubakina, Numerical methods for solving optimization
  problems with differential linear matrix inequalities, Russian Mathematics 64
  (2020) 64--74.

\bibitem{reynolds2020computational}
T.~P. Reynolds, Computational Guidance and Control for Aerospace Systems,
  University of Washington, 2020.

\bibitem{doyle1982analysis}
J.~Doyle, Analysis of feedback systems with structured uncertainties, in: IEE
  Proceedings D (Control Theory and Applications), Vol. 129, IET, 1982, pp.
  242--250.

\bibitem{boyd1994linear}
S.~Boyd, L.~El~Ghaoui, E.~Feron, V.~Balakrishnan, Linear matrix inequalities in
  system and control theory, SIAM, 1994.

\bibitem{accikmecse2003robust}
A.~B. A{\c{c}}ikme{\c{s}}e, M.~Corless, Robust tracking and disturbance
  rejection of bounded rate signals for uncertain/non-linear systems,
  International Journal of Control 76~(11) (2003) 1129--1141.

\bibitem{sontag1995characterizations}
E.~D. Sontag, Y.~Wang, On characterizations of the input-to-state stability
  property, Systems \& Control Letters 24~(5) (1995) 351--359.

\bibitem{accikmecse2002robust}
A.~B. A{\c{c}}{\i}kme{\c{s}}e, M.~Corless, Robust output tracking for
  uncertain/nonlinear systems subject to almost constant disturbances,
  Automatica 38~(11) (2002) 1919--1926.

\bibitem{boyd2004convex}
S.~Boyd, L.~Vandenberghe, Convex optimization, Cambridge {U}niversity {P}ress,
  2004.

\bibitem{oliveira2005stability}
R.~C. Oliveira, P.~L. Peres, Stability of polytopes of matrices via affine
  parameter-dependent {L}yapunov functions: Asymptotically exact {LMI}
  conditions, Linear algebra and its applications 405 (2005) 209--228.

\bibitem{yu2013tube}
S.~Yu, C.~Maier, H.~Chen, F.~Allg{\"o}wer, Tube {MPC} scheme based on robust
  control invariant set with application to {L}ipschitz nonlinear systems,
  Systems \& Control Letters 62~(2) (2013) 194--200.

\bibitem{doi:10.2514/1.G004549}
M.~Szmuk, T.~P. Reynolds, B.~A\c{c}\i{}kme\c{s}e, Successive convexification
  for real-time six-degree-of-freedom powered descent guidance with
  state-triggered constraints, Journal of Guidance, Control, and Dynamics
  43~(8) (2020) 1399--1413.

\end{thebibliography}

\end{document}